\documentclass[11pt]{amsart}
\usepackage{amsmath,amssymb,amsthm}

\newcommand{\sortnoop}[1]{}

\newcommand{\initer}{\setcounter{cste}{\therappel} }
\newcommand{\initre}{\setcounter{rappel}{\thecste} }
\newcounter{rappel}
\newcounter{cste}
\newcommand{\cste}{\addtocounter{cste}{1} c_{\thecste}}
\newcommand{\mcste}{ c_{\thecste}}
\newcommand{\initcste}{\setcounter{cste}{0} }
\newcounter{cstei}
\newcommand{\cstei}{ c_{\thecstei}}
\newcounter{csteii}
\newcommand{\csteii}{ c_{\thecsteii}}
\newcounter{csteiii}
\newcommand{\csteiii}{ c_{\thecsteiii}}
\newcounter{csteij}
\newcommand{\csteij}{ c_{\thecsteij}}
\newcounter{cstej}
\newcommand{\cstej}{ c_{\thecstej}}
\newcounter{csteji}
\newcommand{\csteji}{ c_{\thecsteji}}
\newcounter{cstejii}
\newcommand{\cstejii}{ c_{\thecstejii}}
\newcounter{cstejiii}

\newcounter{cstejij}

\newcounter{cstejj}

\newcounter{cstejji}
\newcommand{\cstejji}{ c_{\thecstejji}}
\newcounter{cstejjii}
\newcommand{\cstejjii}{ c_{\thecstejjii}}
\newcounter{cstejjiii}
\newcommand{\cstejjiii}{ c_{\thecstejjiii}}
\newcounter{cstlem}

\newtheorem{theo}{Theorem}[section]
\newtheorem{cor}[theo]{Corollary}
\newtheorem{prop}[theo]{Proposition}

\newtheorem{lem}[theo]{Lemma}

\newcommand{\delu}{{\rm a}}
\newcommand{\csta}{{\delu}_0}
\newcommand{\cstap}{{C}_0}
\newcommand{\cstaa}{{\rm b}_1}
\newcommand{\cstaaa}{{\rm b'}_1}
\newcommand{\cstb}{{\rm b}_0}

\newcommand{\cstT}{{ \alpha}_0}

\theoremstyle{plain}
\newtheorem{theorem}{Theorem}%[section]
\newtheorem{proposition}[theorem]{Proposition}

\theoremstyle{definition}

\newtheorem*{theorem*}{Theorem}

\newcommand{\reff}[1]{(\ref{#1})}

\newcommand{\findemo}{\hfill $\square$\par
\vspace{0.3cm} 
\noindent}

\newcommand{\cst}{c}

\newcommand{\rw}{{\rm w}}

\newcommand{\cb}{{\mathcal B}}
\newcommand{\cc}{{\mathcal C}l}

\newcommand{\ce}{{\mathcal E}}
\newcommand{\cf}{{\mathcal F}}

\newcommand{\ci}{{\mathcal I}}

\newcommand{\cw}{{\mathcal W}}

\newcommand{\E}{{\mathbb E}}

\newcommand{\N}{{\mathbb N}}
\renewcommand{\P}{{\mathbb P}}

\newcommand{\R}{{\mathbb R}}

\newcommand{\ind}{{\bf 1}}
\newcommand{\range}{{\mathcal R}}

\newcommand{\supp}{{\rm supp}\;}

\newcommand{\val}[1]{\mathop{\left| #1 \right|}\nolimits}
\newcommand{\inv}[1]{\mathop{\frac{1}{ #1}}\nolimits}
\newcommand{\expp}[1]{\mathop {\mathrm{e}^{ #1}}}

\newcommand{\capaf}[1]{\mathop{\mathrm{cap}_f( #1 )}\nolimits}

\begin{document}

 \title[Range of super-Brownian motion]{Some properties of the range of super-Brownian motion}
\author{Jean-Fran\c{c}ois Delmas}
\address{MSRI, 1000 Centennial Drive, Berkeley, CA 94720, U.S.A. \\
and ENPC-CERMICS, 6 av. Blaise Pascal, Champs-sur-Marne, 77455 Marne La Vall\'ee, France.}
\email{delmas@msri.org}
\thanks{The research was done at the \'Ecole Nationale des Ponts et Chauss\'ees and at MSRI, supported by NSF grant DMS-9701755.}

\begin{abstract}

We consider a super-Brownian motion $X$. Its
canonical measures can be studied through
the path-valued process called the Brownian snake. We obtain the
limiting behavior of the volume of the $\varepsilon$-neighborhood for the
range of the Brownian snake, and as a consequence we derive the
analogous result for the range of super-Brownian motion and for the support of the integrated super-Brownian excursion. Then we prove
the support of $X_t$ is capacity-equivalent to 
$[0,1]^2$ in $\R^d$, $d\geq 3$, and the range   of $X$, as well as the support of the integrated super-Brownian excursion are  capacity-equivalent to
$[0,1]^4$ in $\R^d$, $d\geq 5$. 
\end{abstract}

\maketitle 

\section*{Introduction}
Super-Brownian motion, denoted here by  $X=(X_t, t\geq 0)$, is a measure-valued process
in $\R^d$. It can be obtained  as a  limit of branching Brownian particle
systems. We refer to Dynkin \cite{d:ibmp} for such an approximation in a more general
setting. Another way to study     super-Brownian motion, is to use
the path-valued process, called the Brownian snake, which was introduced by
Le Gall \cite{lg:pvmp,lg:pvmppde}. 
Furthermore this approach allows us to study also the integrated super-Brownian excursion (ISE). This process appears naturally when one consider the limit of rescaled lattice trees in high dimension (see Derbez and Slade \cite{ds:ltsbm,ds:smb}).
For every bounded Borel set $A\subset
\R^d$, we denote by $A^\varepsilon=\left\{x\in \R^d; d(x,A)\leq
\varepsilon \right\}$ and  by $\val{A}$   the Lebesgue measure of the set
$A$. 
Recently Tribe \cite{t:rsbm} (see also Perkins \cite{p:smpssbm}) proved a convergence result for
the volume of the $\varepsilon$-neighborhood of the support at time $t>0$, $\supp
X_t$, of  super-Brownian motion in dimension $d\geq 3$.  More precisely,
Tribe showed
that the quantity 
$\varepsilon^{2-d} \val{\left(\supp X_t\right)^\varepsilon\cap A}$
converges a.s. to a
deterministic constant times $\int \ind_A(x) X_t(dx)$. Using results of
Le Gall \cite{lg:hppt} on   hitting
probabilities for the Brownian snake, we give a similar result for the
range of the Brownian snake. We then derive an analogous 
result (theorem \ref{th:rangeX}) for the range of super-Brownian motion after time $t> 0$,
$\range_t(X)$ defined as the closure of $\cup_{s\geq t} \supp X_s$.  
More precisely, we show that there exists a positive constant 
$\cstap$ depending only on $d$ such that for every Borel set $A\subset
\R^d$, $d\geq 4$, for every $t> 0$, we have a.s.
\[
\lim_{\varepsilon\rightarrow 0}
\varphi_d(\varepsilon)\val{\range_t(X)^\varepsilon\cap A}=\cstap
\int_t^\infty  ds\;\int \ind_A(z) X_s(dz),
\]
where $\varphi_4(\varepsilon)=\log (1/\varepsilon)$ and $\varphi_d(\varepsilon)=\varepsilon^{4-d}$ if $d\geq 5$.
We also give a similar result for the support of ISE (corollary \ref{cor:cvpsrangeY}).

Pemantle and Peres \cite{pp:gwt}
defined the notion of capacity-equivalence 
for two random Borel sets, and later  Pemantle and al. \cite{pps:tsbm} showed that the range of 
Brownian motion in $\R^d$, $d\geq 3$, is capacity-equivalent to
$[0,1]^2$. As an application of the previous results , we show
(proposition 
\ref{prop:capaequivX}) that a.s. on $\left\{X_t\neq 0\right\}$, the set $\supp
X_t\subset \R^d$, $d\geq 3$, is capacity-equivalent to $[0,1]^2$, and
that a.s.
the range $\range_t(X)\subset \R^d$ and the support of ISE for $d\geq 5$ are  capacity-equivalent to
$[0,1]^4$.

Let us now describe more precisely the contents of the following
sections. 
In section \ref{sec:setting}, we recall the definition of the
path-valued process $W=(W_s,s\geq 0)$ called the Brownian snake. 
We denote by $\zeta_s$ the lifetime of the path $W_s$. We recall the links between the Brownian snake, super-Brownian motion and ISE.

In section \ref{sec:hitting}, we introduce the main tools concerning the Brownian snake. In particular, we consider $T_{(x,\varepsilon)}$ the hitting time for
the Brownian snake of $\bar B(x,\varepsilon)$, the closed ball with center $x$ and radius
$\varepsilon$:
\[
T_{(x,\varepsilon)}=\inf \left\{ s\geq 0; \exists t\in [0,\zeta_s],
  W_s(t)\in \bar B(x,\varepsilon) \right\}.
\]
The function $u_\varepsilon(x)=\N_0\left[T_{(x,\varepsilon)}<\infty
\right]$, where $\N_0$ is the excursion measure of the Brownian snake away from the trivial path $\rm{0}$, is the maximal nonnegative
solution of $\Delta u =4 u^2$ on $\R^d \backslash B(0,\varepsilon)$ (see also
Dynkin \cite{d:pacnde}).  The study of 
$\val{\range(W)^\varepsilon\cap A}= \int_A dx\;
\ind_{\left\{T_{(x,\varepsilon)}<\infty \right\}}$, where $\range(W)$ is the range of the Brownian snake, relies on the explicit law of the first hitting path  $\left(W_{T_{(x,\varepsilon)}},
  \zeta_{T_{(x,\varepsilon)}}\right)$ under the excursion measure. This law has
been computed by Le Gall \cite{lg:hppt,lg:pccm}. It is closely related
to the law of the process $(x_t^\varepsilon, 0\leq t\leq
\tau^\varepsilon)$, defined as the unique strong solution  of 
\[
dx^\varepsilon_t=d\beta_t+\frac{\nabla
  u_\varepsilon(x^\varepsilon_t-x)}{u_\varepsilon(x^\varepsilon_t-x)}
  dt,\quad \text{for}\quad 0\leq t\leq 
\tau^\varepsilon,
\]
where $\beta$ is a Brownian motion in $\R^d$ started at $\beta_0=0$
and $\tau^\varepsilon=\inf\left\{t\geq 0;
  \val{x^\varepsilon_t-x}=\varepsilon\right\}$. 

In section \ref{sec:rangeX}, we state the main result on the convergence
of the volume of the $\varepsilon$-neighbor\-hood of $\range_t(X)$. The
method of the proof is completely 
different from the one used by Tribe in \cite{t:rsbm}. It  is derived from the convergence of the
volume of the $\varepsilon$-neighborhood of the range of the Brownian snake in
$L^2(\N_{0})$ (proposition \ref{th:erange}).

Section \ref{sec:prprop} is devoted to the proof of the latter
convergence. The proof of the $L^2(\N_{0})$ convergence is somewhat
technical because we need a precise rate of convergence. The derivation
of this estimate relies 
heavily on the explicit law of  $\left(W_{T_{(x,\varepsilon)}},
  \zeta_{T_{(x,\varepsilon)}}\right)$ under $\N_{0}$. It also depends on precise information on the behavior of  the function $u_1$ at infinity. In particular we give the asymptotic expansion of $u_1$ at infinity in the appendix.

In section \ref{sec:cap} we prove the results on capacity-equivalence
for the support and the range of super-Brownian motion and for the support of ISE.
Let $f:[0,\infty )\rightarrow [0,\infty
]$ be a decreasing function. We define
the energy of a Radon measure $\nu$ on $\R^d$ with respect to the
kernel $f$ by: $ \ci_f(\nu)=\iint f(\val{x-y}) \nu(dx)\nu(dy)$,
and the  capacity of a set $\Lambda\subset \R^d$ by 
$
\capaf{\Lambda} =\left[\inf_{\nu(\Lambda)=1} \ci_f(\nu)\right] ^{-1}
$. 
Following the terminology introduced  in
\cite{pp:gwt}, we say
that two sets $\Lambda_1$ and
$\Lambda_2$ are capacity-equivalent 
if there exist two positive constants $c$ and $C$ such that for every kernel $f$, we have
\[
c\capaf{\Lambda_1} \leq \capaf{\Lambda_2}\leq C\capaf{\Lambda_1} .
\] 
Proposition \ref{prop:capaequivX} states that a.s. the set $\supp
X_t\subset \R^d$, $d\geq 3$, is capacity-equivalent to $[0,1]^2$, and
that a.s.
the range $\range_t(X)\subset \R^d$, as well as the support of ISE for $d\geq 5$ are  capacity-equivalent to
$[0,1]^4$. The proof follows the method of
\cite{pps:tsbm}. 

\section{Preliminaries  on the Brownian snake and super-Brownian motion}
\label{sec:setting}
We first introduce some notation. We denote by $(M_f,{\mathcal M}_f)$ the
space of all finite
measures on $\R^d$, endowed with the  topology of  weak 
convergence.  We denote by  $\cb_{b+} (\R^p)$, respectively $\cb_{b+} (\R^+\times \R^p)$,
the set of all real bounded nonnegative measurable functions defined on $\R^p$,
respectively on $\R^+\times \R^p$.  We also 
denote by $\cb (\R^p)$  the Borel $\sigma$-field  on $\R^p$. For $A\in \cb
(\R^p)$, let $\cc (A)=\bar A$ be the  closure of 
$A$. For every measure $\nu\in M_f$,
and 
$f\in \cb_{b+}(\R^d)$, we shall write $\int f(y)\nu(dy)
=(\nu,f)$. We also denote by $\supp \nu$ the closed support of the
 measure $\nu$. If $S$ is a Polish space, we denote by $C(I, S)$ the set of
all continuous functions from $I\subset \R$ into $S$.

\subsection{The Brownian snake}

\label{sec:snake}
We recall some facts about the Brownian snake, a path-valued Markov
process introduced by Le Gall \cite{lg:pvmp,lg:pvmppde}. A stopped path
is a continuous function $\rw :[0,\zeta]\rightarrow \R^d$, where
$\zeta=\zeta_{(\rw)}$ is called the lifetime of the path.
We shall denote by $\hat \rw $ the end point $\rw(\zeta)$. Let $\cw$ be
the space of all stopped  paths in $\R^d$. When equipped with the metric
\[
d(\rw,\rw')=\val{\zeta_{(\rw)} -\zeta_{(\rw')}}+\sup_{s\geq
  0}\val{\rw(s\wedge\zeta_{(\rw)})-\rw'(s\wedge\zeta_{(\rw')})},
\]
the space $\cw$ is a Polish space.

Let $\rw \in \cw$ and $a,b\geq 0$, such that $a\leq
b\wedge\zeta_{(\rw)}$. There exists a unique probability measure on $\cw$
denoted by $Q_{a,b}^\rw(d\rw')$ such that:
\begin{itemize}
   \item[(i)] $\zeta_{(\rw')}=b$, $Q_{a,b}^\rw(d\rw')$-a.s.
   \item[(ii)] $\rw'(t)=\rw(t)$ for every $0\leq t\leq a$,
     $Q_{a,b}^\rw(d\rw')$-a.s.
   \item[(iii)] The law of $\left(\rw'(t+a), 0\leq t\leq b-a
     \right)$ under $Q_{a,b}^\rw(d\rw')$ is the law of Brownian
     motion in $\R^d$ started at $\rw (a)$ and stopped at time $b-a$.
\end{itemize}
We shall also consider $Q_{a,b}^\rw(d\rw')$ as a probability on the
space %of continuous functions 
$C([0,b],\R^d)$.  We set  $\cw_x=\left\{\rw \in \cw; \rw(0)=x
\right\}$ for $x\in
\R^d$. Let $\rw \in \cw_x$. We restate theorem 1.1 from
\cite{lg:pvmp}:
\begin{theo}[Le Gall] 
\label{th:serpentlg}
  There exists a continuous strong
  Markov process with values in
  $\cw_x$, $W=(W_s,s\geq 0)$,  whose law is
  characterized by the following two properties. 
\begin{itemize}
   \item[(i)] The lifetime process
     $\zeta=\left(\zeta_s=\zeta_{(W_s)},s\geq 0\right)$ is a reflecting 
     Brownian motion in $\R^+$.% which starts at $\zeta_0=\zeta_{(\rw)}$.
   \item[(ii)] Conditionally given $\left(\zeta_s,s\geq 0\right)$, the
     process $\left(W_s,s\geq 0\right)$ 
     is a time-inhomogeneous continuous
     Markov process, % which starts at
     %$W_0=\rw$ and 
     whose transition kernel between times
     $s$ and $s'\geq s$ is 
\[
P_{s,s'}(\rw,d\rw')=Q_{m(s,s'),\zeta_{s'}}^\rw(d\rw'),
\]
where $m(s,s'):=\inf_{r\in [s,s']}\zeta_r$.
\end{itemize}
\end{theo}
From now on we shall consider the canonical realization of the process
$W$ defined on the space $C(\R^+,\cw_x)$. The law of $W$ started at
$\rw$ is denoted by $\ce_{\rw}$.
We will use the following
consequence of (ii): outside a $\ce_\rw$-negligible set, for every
$s'>s$, one has 
$W_s(t)=W_{s'}(t)$ for every $t\in [0,m(s,s')]$. We
shall write $\ce^*_\rw$ for the law of the process $W$ killed when its
lifetime reaches zero. The  distribution of $W$ under $\ce^*_\rw$ can be characterized as in
theorem \ref{th:serpentlg}, except that its lifetime process is
distributed as a linear Brownian motion killed at its first hitting time of
$\{0\}$. The state space for $\left(W,\ce^*_\rw\right)$ is the space
$\cw_x^*=\cw_x\cup \partial$, where $\partial$ is a cemetery point. The trivial
path $\textbf{x}$ such that $\zeta_{(\textbf{x})}=0$,
$\textbf{x} (0)=x$ is clearly a regular point for the process
$\left(W, \ce_{\rw}\right)$. Following \cite{b:emp} chapter 3, we can consider the
excursion measure, $\N_x$, outside $\left\{\textbf{x}\right\}$. 
The
distribution of $W$ under $\N_x$ can be characterized as in theorem
\ref{th:serpentlg}, except that now the lifetime process $\zeta$ is
distributed according to It\^o measure of positive excursions of linear
Brownian motion. We  normalize $\N_x$ so that, for every
$\varepsilon>0$,
\[
\N_x\left[\sup_{s\geq 0}\zeta_s>\varepsilon\right]=\inv{2\varepsilon}.
\]
The Brownian snake enjoys a scaling property: if $\lambda>0$, the law of the process  
$W^{(\lambda)}_s(t)= \lambda^{-1} W_{\lambda^4 s} (\lambda^2 t)$ under
$\N_x$ is $\lambda^{-2} \N_{\lambda^{-1}x}$.

We recall the strong Markov property for the snake under $\N_x$ (see
\cite{lg:pvmppde}). Let $T$ be a stopping time of the natural filtration
$\cf^W$ of the process $W$. Assume $T>0$ $\N_x$-a.e., and let $F$, $H$ nonnegative measurable
functionals on $C(\R^+,\cw_x^*)$ such that $F$ is $\cf_T^W$ measurable. Then
if $\theta$ denotes the usual shift operator, we have 
\[
\N_x\left[T<\infty ; F\:\cdot\: H\circ \theta_T\right] =
\N_x\left[T<\infty ;F\:\cdot\: \ce_{W_T}^*\left[H\right]\right]. 
\]
Let $\sigma=\inf\left\{s>0;\zeta_s=0\right\}$ denote  the duration of the
excursion of $\zeta$ under $\N_x$. The range $\range =\range (W)$ of $W$ is defined under
$\N_x$ by
\[
\range=\left\{W_s(t); 0\leq t\leq \zeta_s, 0\leq  s \leq  \sigma\right\} =
\left\{\hat W_s; 0\leq  s \leq  \sigma \right\}.
\]

For every nonnegative measurable function $F$ on $\cw_x^*$, we have
\begin{equation*}
%\label{eq:NGds}
\N_x\left[\int_0^\sigma F(W_s,\zeta_s) ds \right] = \int_0^\infty
\E_x\left[F(\beta_{[0,t]},t)\right] dt,
\end{equation*}
where $\beta_{[0,t]}$ is under $\P_x$ the restriction to $[0,t]$ of a Brownian motion in $\R^d$ started at
$\beta_0=x$. Now consider under $\N_x$ the continuous version
$\left(l^t_s,t> 0,s\geq 0\right)$ of the local time of $\zeta$ at
level $t$ and time $s$. We define a  measure valued process $Y$ on
$\R^d$ by setting  for every $t>0$, for every $\varphi\in \cb_{b+}(\R^d)$, 
\[
(Y_t,\varphi)=\int_0^\sigma dl_s^t \; \varphi(\hat W_s).
\]
We shall sometimes write $Y_t(W)$ to recall that $Y_t$ is a function of
the Brownian snake. From the joint continuity of the local time and the continuity of the map $s\mapsto \hat W_s$, we get that $\N_x$-a.e., the process 
$Y$ is continuous on $(0,\infty )$ for the Prohorov 
distance on $M_f$. Let $\varphi\in
\cb_{b+}(\R^d)$. We define  on $\R^+\times \R^d$ the function
$v(t,x)=\N_x\left[1-\exp{-(Y_t,\varphi)}\right]$, if $t>0$, and $v(0,x)= 
\varphi(x)$. We will write $v(t)$ for the function $v(t,\cdot)$. We recall that the
function $v$   is the unique  nonnegative measurable 
solution of the integral functional equation
\begin{equation}
\label{eq:integrale}
v(t)+2\int_0^t ds\; P_s\left[v(t-s)^2\right]=J(t)\qquad t\geq 0,
\end{equation}
where $J(t,x) =P_t [\varphi](x)$, and $(P_t,t\geq 0)$ is the Brownian semi-group
in $\R^d$. A few other  remarks on the solution of
(\ref{eq:integrale}) are presented in section \ref{sec:annexemo} below. 

\subsection{Super-Brownian motion and ISE}

Let us now recall the 
definition of  super-Brownian motion  %(see \cite{d:bpss,d:ibmp}) 
and its connection with the Brownian snake. %(see \cite{lg:pvmp,lg:pvmppde})
The second part of the next 
theorem is lemma 4.1 from \cite{d:bpss}. 
Let $\nu\in M_f$. 
\begin{theo}
   %Let $\nu\in M_f$. 
   There exists a continuous strong Markov process $X=\left(X_s,s\geq
   0\right)$ defined on the canonical space $C(\R^+, M_f)$,
   whose law is characterized by the two following properties under $ \P^X_\nu$.
\begin{itemize}
   \item[(i)] $X_0=\nu$, $\P^X_\nu$-a.s.
   \item[(ii)] For every $\varphi\in \cb_{b+}(\R^d)$, $t\geq s>0$, we have
\begin{equation*}
%\label{eq:LaplaceX}
\E_\nu^X\left[\exp{\left[-(X_t,\varphi)\right]}\mid \sigma(X_u,0\leq
    u\leq s)\right] = \exp{\left[-(X_s,v(t-s))\right]} ,
\end{equation*}
    where the function $v$ is the unique nonnegative solution of
    (\ref{eq:integrale}) with $J(t)=P_t[\varphi]$.
\end{itemize}
Furthermore, for every integer $m\geq 1$, $t_m>\cdots>t_1\geq 0$,
$\varphi_1,\ldots,\varphi_m \in \cb_{b+}(\R^d)$, we have
\begin{equation}
\label{eq:fLaplaceX}
\E_\nu^X\left[\exp{\left[-\sum_{\left\{i;t_i\leq
t\right\}}(X_{t-t_i},\varphi_i)\right]}\right] =
\exp{\left[-(\nu,v(t))\right]} ,
\end{equation}
where $v$ is the unique nonnegative solution to the integral equation
(\ref{eq:integrale}) with right-hand side  $J(t)=\sum_{\left\{i;t_i\leq
t\right\}}P_{t-t_i}[\varphi_i]$.
\end{theo}
\begin{theo}[Le Gall \cite{lg:pvmp,lg:pvmppde}]
\label{th:XY}
   Let $\sum_{i\in I} \delta_{W^i}$ be a Poisson measure on
   $C(\R^+,\cw)$ with intensity $\int \nu(dx) \N_x[\cdot]$, then the process
   $Z$ defined by $Z_0=\nu$ and  $Z_t=\sum_{i\in I}Y_t(W^i)$ if $ t>0$,
   is distributed according to $\P_\nu^X$.
\end{theo}
We deduce from the normalization of $\N_x$ that $\N_x\left[Y_t\neq
0\right]=1/2t<\infty $. This implies that for every $t>0$, there
is only a finite number of indices $i\in I$ such that the process $(Y_s(W^i), 
s\geq t )$ is nonzero.

We now recall the connection between ISE and Brownian snake.
There exists a unique collection $\left(\N^{(r)}_0, r>0\right)$ of probability measure on $C(\R^+,\cw_0^*)$ such that:
\begin{enumerate}
   \item For every $r>0$, $\N^{(r)}_0[\sigma=r]=1$.
   \item For every $\lambda>0$, $r>0$, $F$, nonnegative measurable functional on $C(\R^+,\cw_0^*)$,
\[
\N^{(r)}_0\left[F(W^{(\lambda)})\right]=\N^{(\lambda^{-4} r)}_0\left[F(W)\right].
\]
   \item For every nonnegative measurable functional $F$ on $C(\R^+,\cw_0^*)$,
\begin{equation}
   \label{eq:Nnormal}
\N_0[F]=\inv{\sqrt{2\pi}} \int_0^\infty  dr\; r^{-3/2} \N^{(r)}_0[F].
\end{equation}
\end{enumerate}
The measurability of the mapping $r\mapsto \N^{(r)}_0[F]$ follows from the scaling property 2. Under $\N^{(1)}_0$, the distribution of $W$ is characterized as in theorem \ref{th:serpentlg}, except that the lifetime process is distributed according to the normalized It\^o measure. The law of the ISE is the law of the continuous tree associated to $\sqrt{2}W$, under $\N^{(1)}_0$ (see corollary 4 in \cite{lg:urtbe} and  \cite{a:tbm}). In particular the law of the support of ISE is the law of $\sqrt{2}\range$ under $\N^{(1)}_0$, where we set $\lambda A=\{x; \lambda^{-1} x\in A\}$.

\subsection{Hitting probabilities for the Brownian snake}
\label{sec:hitting}

We now recall a few  results from
\cite{lg:hppt}. 
Let $\rw \in \cw\cup C(\R^+,\R^d)$,
we introduce the first hitting time of $A\in \cb(\R^d)$:
\[
\tau_{A}(\rw)=\inf \left\{t\geq 0; \rw(t)\in A \right\},
\]
with the usual convention $\inf \emptyset =\infty $.
We   omit $\rw$ 
when there is no risk of confusion.   Consider the Brownian snake $W$, and set 
\[
T_{(y,\varepsilon)}=\inf \left\{s\geq 0; \exists t\in [0,\zeta_s],
  W_s(t)\in \bar B(y,\varepsilon) \right\},
\]
where $B(y,\varepsilon)$ is the
open ball in $\R^d$ centered at 
$y$ with radius $\varepsilon> 0$, and $\bar B(y,\varepsilon)$ its
closure. We know from \cite{lg:pvmppde} that the function defined on
$\R^d\backslash \bar B(0,\varepsilon)$,
\[
u_\varepsilon(y):=\N_0\left[T_{(y,\varepsilon)}<\infty
\right]=\N_0\left[\range \cap \bar B{(y,\varepsilon)}\neq \emptyset
\right]=\N_{-y}\left[\range \cap \bar B{(0,\varepsilon)}\neq \emptyset
\right] , 
\]
is the maximal nonnegative solution on $\R^d\backslash \bar
B(0,\varepsilon)$ of 
\[
\Delta u= 4u^2.
\]
This result was first proved in a more general setting by Dynkin
\cite{d:pacnde} in terms of superprocesses. The function $u_\varepsilon$
is strictly positive on $\R^d\backslash \bar B(0,\varepsilon)$. For
every $y_0\in \partial B(0,\varepsilon)$, we have 
\[
\lim_{y\in \bar B(0,\varepsilon)^c; y\rightarrow y_0} u_\varepsilon(y)=\infty .
\]
Scaling and symmetry arguments show that for every $y\in \R^d\backslash\bar B(0,\varepsilon)$,
\begin{equation}
\label{eq:uscale}   
u_\varepsilon(y)=\varepsilon^{-2}
u_1\left(\frac{\val{y}}{\varepsilon}\right),
\end{equation}
where the function $u_1(r)$, $r\in (1,\infty )$ is the maximal nonnegative solution on
$(1,\infty )$ of
\[
u''_1(r)+\frac{d-1}{r} u'_1(r) =4 u^2_1(r).
\]
It is easy to see that the function $u_1$ is decreasing. 
In section \ref{sec:annexeu} we give the asymptotic expansion of $u_1$ at infinity.

We give the following result on the probability of the event
$\left\{T_{(y,\varepsilon)}<\infty \right\}$ (see lemma 2.1 of
\cite{lg:hppt}). Assume $x_0\not \in \bar B(y,\varepsilon)$. Then $\N_{x_0}$-a.e. for every
$T \geq 0$, we have
\begin{align}
\label{eq:hitting}
\ce^*_{W_T} \left[T_{(y,\varepsilon)}<\infty \right]
&=2\int_0^{\zeta_T\wedge \tau_{B(y,\varepsilon)}(W_T)}
dt\;u_\varepsilon(W_{T}(t)-y) ) 
\expp{\left[-2\int_0^t u_\varepsilon(W_{T}(s)-y) ds \right]}\\
\nonumber
&=1-
\exp{\left[-2\int_0^{\zeta_T\wedge \tau_{B(y,\varepsilon)}(W_T)} u_\varepsilon(W_{T}(s)-y) ds \right]}.
\end{align}

Let $x_0,x\in \R^d$. We will now describe the law of
$W_{T_{(x,\varepsilon)}}$ under $\N_{x_0}[\cdot\mid T_{(x,\varepsilon)}<\infty ]$.
 First of all we denote by
$\beta$ a Brownian motion in $\R^d$ started at $x_0$ under
$\P_{x_0}$. Assume $x_0\not \in \bar B(x,\varepsilon )$. Corollary 2.3 from
\cite{lg:hppt} ensures that there exists 
$\P_{x_0}$-a.s. a
unique continuous process $x^\varepsilon= (x^\varepsilon_t, 0\leq
t\leq \tau^\varepsilon)$ 
taking values in $\R^d$ such that for every $\eta\in(0, \val{x-x_0}
-\varepsilon)$, for every $t\leq  \tau_\eta^\varepsilon=\inf \left\{s\geq 0; \val{x^\varepsilon_s-x}
  \leq \varepsilon+ \eta \right\}$,
\[
x^\varepsilon_t=\beta_t +
\int_0^{t} \;\frac{\nabla
u_\varepsilon (x^\varepsilon_s-x)}{u_\varepsilon (x^\varepsilon_s-x)} ds,
\]
furthermore,  $\P_{x_0}$-a.s. $\tau^\varepsilon=\lim_{\eta\rightarrow 0}
\tau_\eta^\varepsilon<\infty $ and
$\val{x^\varepsilon_{\tau^\varepsilon} -x}= \varepsilon$. 
We also recall that thanks  to Girsanov's theorem, we have
for every nonnegative measurable function $F$ on $C([0,t], \R^d)$
\begin{multline*}
\E_{x_0} \left[\tau^\varepsilon>t; F\left(x^\varepsilon_{[0,t]}\right)\right] \\
=\E_{x_0}\left[\tau_{B(x,\varepsilon)}(\beta)>t;F\left(\beta_{[0,t]}
\right)\frac{u_\varepsilon(\beta_t-x)}{u_\varepsilon(x_0-x)}\exp{\left[-2\int_0^t
u_\varepsilon(\beta_s-x) ds  \right]}
\right],
\end{multline*}
where $x^\varepsilon_{[0,t]}$ and $\beta_{[0,t]}$ are the
restriction of $x^\varepsilon$ and $\beta$ to $[0,t]$.
The law of $x^\varepsilon$ under $\P_{x_0}$ can be interpreted as a
probability measure on $\cw_{x_0}^*$. Consider the closed set 
\[
A=\left\{\rw\in \cw_{x_0}^*; \tau_{\bar B(x,\varepsilon)}(\rw)<\infty \right\}.
\]
It has been proved in \cite{lg:hppt} that its capacitary measure with
respect to the Brownian snake with initial point  $x_0$ is 
exactly $u_\varepsilon(x_0-x)$ times the law of $x^\varepsilon$ under
$\P_{x_0}$. 
It is not hard to check however that the normalized capacitary measure can be
interpreted as the hitting distribution under $\N_{x_0}$ (cf
\cite{lg:pccm}, this is proved in a way similar to the classical
interpretation of the capacitary measure as a last exit distribution,
see e.g. Port and Stone \cite{ps:bm}). 
Thus we deduce that for every
nonnegative measurable function $F$ on $\cw_{x_0}^*$, we have
\begin{equation*}
%\label{eq:Nnu}
\N_{x_0}\left[T_{(x,\varepsilon)}<\infty ; F(W_{T_{(x,\varepsilon)}},\zeta_{T_{(x,\varepsilon)}})\right]
=u_\varepsilon(x_0-x)\E_{x_0} \left[F(x^\varepsilon,\tau^\varepsilon)\right].
\end{equation*}
Hence for every $t\geq 0$, and every nonnegative measurable function $F$
on $C([0,t], \R^d)$, we have
\begin{multline}
\N_{x_0}\left[T_{(x,\varepsilon)}<\infty ; 
\zeta_{T_{(x,\varepsilon)}}>t; 
F\left((W_{T_{(x,\varepsilon)}}(s),s\in [0,t])\right)\right]\\
\label{eq:Nbeta}
=\E_{x_0}\left[\tau_{B(x,\varepsilon)}>t; F\left(\beta_{[0,t]}\right){u_\varepsilon(\beta_t-x)}\exp{\left[-2\int_0^t
u_\varepsilon(\beta_s-x) ds  \right]}
\right].
\end{multline}
Finally we shall use the
following inequality, that can be derived from the 
Feynman-Kac  formula (use the fact that $u_\varepsilon$ solves $\Delta
u=4u_\varepsilon u$)
\begin{equation}
\label{eq:fk}
\hspace{-0.5cm}u_\varepsilon(x)\geq  2\E_0 \left[ \int_0^{\tau_{B(x,\varepsilon)} }dt\;\;
    {u_\varepsilon(\beta_t-x)}^2\exp{\left[-4\int_0^t 
u_\varepsilon(\beta_s-x) ds  \right]}\right].
\end{equation}
There is in fact equality in \reff{eq:fk} (see the  remark on  page 293 of \cite{lg:hppt}).

\section{A  property of the range of  super-Brownian motion}
\label{sec:rangeX}
\initcste

\noindent
For $A\in \cb(\R^d)$,
$\varepsilon>0$, we set  $ A^\varepsilon:=\left\{x\in \R^d; d(x,A)\leq \varepsilon \right\}$,
with $d(x,A)=\inf \left\{ \val{x-y}; y\in A\right\}$.
We will write $\val{A}$ for the Lebesgue measure of $A$. We also set 
\[
\cstap=\csta 2\pi^{d/2} \Gamma([d-2]/2)^{-1},
\]
where the
constant $\csta$ is defined in lemma
\ref{lem:a0} (see also the remark below the lemma). 
We set $\range_t(X)=\cc \left(\bigcup _{s\geq
  t} \supp X_s\right)$. Let $\varphi_d(\varepsilon)=\varepsilon^{4-d}$ if $d\geq 5$ and $\varphi_4(\varepsilon)=\log (1/\varepsilon)$ for $\varepsilon>0$.
\begin{theo}
\label{th:rangeX} 
Let $\nu\in M_f$. For every  Borel set $A\subset \R^d$, $d\geq 4$,
for every $t>0$, $\P_\nu^X$-a.s.
\begin{equation}
   \label{eq:rangeX}
\lim_{\varepsilon\rightarrow 0}
\varphi_d(\varepsilon)\val{\range_t(X)^\varepsilon\cap A}=\cstap
\int_t^\infty  ds\;(X_s,\ind_A).
\end{equation}
If there exists $\rho<4$ such that
$\lim_{\varepsilon\rightarrow 0} \varepsilon^{\rho-d}\val{(\supp
  \nu)^\varepsilon}=0$ then \reff{eq:rangeX} holds  with $t=0$.
\end{theo}
Let $K$ a compact subset of $\R^d$. We consider the measure $\phi(K)$ defined by $\phi(K)(A)=\val{K\cap A}$. Since the set $\range_t(X)$ is compact for $t>0$, the theorem implies that a.s. the sequence of measures $(\varphi_d(\varepsilon) \phi(\range_t(X)^\varepsilon), \varepsilon>0)$ converges weakly to $\cstap
\int_t^\infty  ds\;(X_s,\ind_A)$. 

Let us recall the main theorem of 
\cite{t:rsbm} (see also \cite{p:smpssbm}).
\begin{theo}[Tribe]
\label{th:tribe}
Let $A$ a bounded Borel set in $\R^d$, $d\geq 3$. Fix $t>0$ and $\nu\in
M_f$. Then there exists a positive constant $\cstT$ depending only on
$d$ such that 
\[
\lim_{\varepsilon\rightarrow 0} \varepsilon^{2-d}\val{\left(\supp
X_t\right)^\varepsilon\cap A} =\cstT (X_t,\ind_A),
\]
where the convergence holds $\P^X_\nu$-a.s. and in $L^2(\P^X_\nu)$.
\end{theo}
We
shall deduce theorem \ref{th:rangeX} from the next proposition on the range of the Brownian
snake, whose proof will be given in  the next section. For $\theta\in (0,1)$, we set $h_{d,\theta}(\varepsilon)=\varepsilon^{1-\theta}$ if $d\geq 5$ and $h_{4,\theta}(\varepsilon)=\log (1/\varepsilon)^{-1/\theta}$ for $\varepsilon\in (0,1)$. For short we will write $h_d$ for $h_{d,\theta}$.
\begin{prop}
\label{th:erange}
Let $d\geq 4$. For every $\theta\in (0,1/d)$ and every $R_0>0$, there
exists a constant $\kappa=\kappa(\theta)>0$ and  $\varepsilon_0>0$
such that for every $\varepsilon\in (0,\varepsilon_0]$, for every
$x_0$ with $\val{x_0}\leq R_0$,  and every Borel set $A\subset \bar B(0,R_0)$, we have 
\[
\val{\N_{x_0}\left[\varphi_d(\varepsilon)\val{\range(W)^\varepsilon\cap A\cap \bar
  B(x_0,h_d(\varepsilon))^c }-\cstap \int_0^\infty
  ds\;(Y_s,\ind_A)\right]}\leq h_d(\varepsilon)^{\kappa/2},
\]
and
\[
\N_{x_0}\left[\left[\varphi_d(\varepsilon)     \val{\range(W)^\varepsilon\cap A\cap \bar
  B(x_0,h_d(\varepsilon))^c }-\cstap \int_0^\infty
  ds\;(Y_s,\ind_A)\right]^2\right]\leq h_d(\varepsilon)^\kappa.
\]
\end{prop}
\noindent
\textbf{Remark}. We have trivially $B(x_0,\varepsilon)\subset
\range(W)^\varepsilon$, $\N_{x_0}$-a.e. Since $\N_{x_0}$ is an infinite
measure, $\N_{x_0} \left[ \val{\range(W)^\varepsilon\cap
    B(x_0,\delta)}\right]=\infty $ for every
$\varepsilon,\delta>0$. This is the reason why we consider $A\cap \bar
  B(x_0,h_d(\varepsilon))^c $ rather than $A$ in the previous
  proposition.

We first give a consequence of this proposition.
\begin{cor}
\label{cor:cvpsrangeY}
   Let $d\geq 4$. For every Borel set $A\subset \R^d$, $\N_{x_0}$-a.e., we have 
\[
\lim_{\varepsilon\rightarrow 0} \varphi_d(\varepsilon)\val{\range(W)^\varepsilon \cap
  A}=\cstap \int_0^\infty ds\; (Y_s,\ind_A).
\]
The results holds $\N^{(1)}_0$-a.s. if $\val {\partial A}=0$.
\end{cor}
\noindent
\textbf{Proof} of corollary \ref{cor:cvpsrangeY}. Since
$\N_{x_0}$-a.e. the range $\range(W)$ is bounded, we only need  to consider
a bounded Borel set $A$. Let $\kappa>0$ be fixed as in proposition \ref{th:erange}.
Let $\varepsilon_n$ such that $h_d(\varepsilon_n)=n^{-2/\kappa}$ for $n\geq 1$. 
Using the Borel-Cantelli lemma and the second
upper bound of proposition \ref{th:erange}, we get that
 the sequence
$\left(\varphi_d(\varepsilon_n)\val{\range(W)^{\varepsilon_n} \cap A},n\geq
 1\right)$ converges $\N_{x_0}$-a.e. to $\cstap
\int_0^\infty  ds\; (Y_s,\ind_A)$. But for $\varepsilon'\leq \varepsilon$,
since $\range(W)^{\varepsilon'} \subset \range(W)^\varepsilon$, we have
\[
\varphi_d({\varepsilon'}) \val{\range(W)^{\varepsilon'}\cap A} \leq
\varphi_d({\varepsilon}) \val{\range(W)^\varepsilon \cap
  A}\varphi_d({\varepsilon'})/\varphi_d({\varepsilon}).
\]
A monotonicity argument using the fact that
$\varphi_d({\varepsilon_{n+1}})/\varphi_d({\varepsilon_n})$ converges to $1$, completes the proof of the first part.

The above result implies that $\N_0$-a.e. the sequence of measures $(\varphi_d(\varepsilon) \phi(\range(W)^\varepsilon),\varepsilon>0)$ converges weakly to $\cstap \int_0^\infty ds\; Y_s$.
Using \reff{eq:Nnormal} we see this convergence also holds $dr$-a.e. $\N^{(r)}_0$-a.s. By the scaling property the Brownian snake and  the family $(\N^{(r)}_0, r>0)$, we get this convergence holds $\N^{(1)}_0$-a.s. Thus we have for every Borel set $A\subset \R^d$, $\N^{(1)}_0$-a.s.
\begin{multline*}
   \cstap \int_0^\infty ds\; (Y_s,\ind_{\text{Int}(A)})\leq \liminf_{\varepsilon\rightarrow 0} \varphi_d(\varepsilon) \val{\range(W)^\varepsilon\cap A}\\
\leq \limsup_{\varepsilon\rightarrow 0} \varphi_d(\varepsilon) \val{\range(W)^\varepsilon\cap A} \leq  \cstap \int_0^\infty ds\; (Y_s,\ind_{\bar A}),
\end{multline*}
where Int$(A)$ denotes the interior of $A$. To prove the second part of the corollary we just need to check that if $\val{\partial A}=0$ then $\int_0^\infty ds\; (Y_s,\ind_{\text{Int}(A)})=\int_0^\infty ds\; (Y_s,\ind_{\bar A})$. It is enough to prove that $\val{A}=0$ implies $\int_0^\infty ds\; (Y_s,\ind_{A})=0$ $\N^{(1)}_0$-a.s. Conditioning by the lifetime process, we get
\[
\N^{(1)}_0\left[\int_0^\infty ds\; (Y_s,\ind_{ A})\right]=\N^{(1)}_0\left[\int_0^1 dt\; \ind_{A}(\hat W_t)\right]=\int_0^1 dt \;\N^{(1)}_0\left[P_{\zeta_t}[\ind_A](0)\right].
\]
This is equal to zero if $\val{A}=0$. This ends the proof of the second part of the corollary.
\findemo
As a byproduct of the proof we get that  $\N_{x_0}$-a.e. and $\N^{(1)}_0$-a.s. the sequence of measures $(\varphi_d(\varepsilon) \phi(\range(W)^\varepsilon),\varepsilon>0)$ converges weakly to $\cstap \int_0^\infty ds\; Y_s$.

We first state some straightforward consequences of \reff {eq:uscale} and lemma \ref{lem:a0}. We say that $\varepsilon_0>0$ satisfies the condition (C) if $\varepsilon_0^{-\theta}\geq 4/3$ if $d\geq 5$ or  $\log (1/\varepsilon_0)\geq  4 \log (2/\theta)/\theta$ if $d=4$.  For $\theta\in (0,1/4)$ this implies that for $\varepsilon\in (0,\varepsilon_0)$, $h_4(\varepsilon)/\varepsilon\geq 4/3$ and that 
\begin{equation}
   \label{eq:lolo}
\log (\log (1/\varepsilon))/[\theta \log (1/\varepsilon)]\leq 1/2.
\end{equation}
For $d\geq 4$, $\theta\in (0,1/d)$, there  exists a constant $\cstaa$ such that for every $\varepsilon$ satisfying (C), $x\not \in B(0,h_d(\varepsilon))$ we have 
\begin{align}
  \label{eq:umajo}  
u_\varepsilon(x)&\leq \cstb \varphi_d(\varepsilon)^{-1} \val{x}^{2-d},\\
  \label{eq:udel}
u_\varepsilon(x)&\leq  \varphi_d(\varepsilon)^{-1} \val{x}^{2-d}\left[\csta+\cstaa h_d(\varepsilon)^{\theta/2}\right]. 
\end{align}
For $\val{x}>\varepsilon$, we have 
\begin{align}
\label{eq:umino}
   &u_\varepsilon(x)\geq \csta \varphi_d(\varepsilon)^{-1} \val{x}^{2-d} \quad \text{if $d\geq 5$},\\ 
\label{eq:umino4}
& u_\varepsilon(x)\geq \csta \varphi_4(\varepsilon)^{-1} \val{x}^{-2} \left[1+\log (2\val{x})/\log (1/\varepsilon)\right]^{-1} \quad \text{if $d=4$}.
\end{align}
We will also often use the following inequality for $\varepsilon$ satisfying (C): $\varphi_d(\varepsilon)h_d(\varepsilon)^d\leq h_d(\varepsilon)^3$.

\noindent
\textbf{Proof} of theorem \ref{th:rangeX}. Recall that for every
$t>0$, $\P^X_\nu$
a.s. the set  $\range_t(X)$ is bounded.  Thus we only need to consider a  bounded
Borel set $A$. 
Thanks to
the Markov property of $X$ at time $t$ and  theorem \ref{th:tribe}
it is clearly enough to prove the second part of theorem \ref{th:rangeX}. Let $\nu\in M_f$
and $\rho<4$ such that $\lim_{\varepsilon\rightarrow 0} \varepsilon^{\rho-d}\val{(\supp
  \nu)^\varepsilon}=0$. For short we write a.s. for $\P^X_\nu$-a.s.

First step. Recall   we
can write for every $t>0$, $X_t=\sum_{i\in I} Y_t(W^i)$, where
$\sum_{i\in I} \delta_{W^i}$ is a Poisson measure on
$C(\R^+, \cw)$ with intensity measure $\int \nu(dx) \N_x[\cdot]$.
We let $x_0^i$ denote the starting point of the Brownian snake $W^i$
(i.e. $x_0^i= W^i_0(0)$). Notice that a.s. for every $i\in I$, $x^i_0\in
\supp \nu$, which is bounded thanks to the hypothesis on $\supp \nu$.
Fix $\theta\in (0,1/d)$ such that $d-\rho\geq
(d-4)/(1-\theta)$ (and $\theta<4-\rho$ if $d=4$). Fix  $R_0$ such that $ \supp \nu\subset
B(0, R_0)$.  Let $\kappa$ and
$\varepsilon_0<1$ fixed as in proposition \ref{th:erange}. We notice that for every bounded Borel
set $A\subset B(0, R_0)$, 
\[
\varphi_d(\varepsilon)\val{\range_0(X)^\varepsilon\cap A}\leq \sum_{i\in I}
V_\varepsilon(W^i) + \varphi_d(\varepsilon)\val{A\cap (\supp \nu)^{h_d(\varepsilon)}},
\]
where 
\[
V_\varepsilon(W^i)=\varphi_d(\varepsilon)\val{\range(W^i)^\varepsilon\cap A\cap \bar
  B(x_0^i,h_d(\varepsilon))^c }.
\]
We set $V_0(W^i)=\cstap \int_0^\infty
  ds\;(Y_s(W^i),\ind_A)$.
We use the second moment formula  for a Poisson measure to get:
\begin{multline*}
\E^X_\nu\left[\left[\sum_{i\in I}
V_\varepsilon(W^i)-\sum_{i\in I}
V_0(W^i)\right]^2\right]
=\int \nu(dx) \N_x \left[\left[
V_\varepsilon(W)-V_0(W)\right]^2\right]\\
+\left[\int \nu(dx) \N_x \left[
V_\varepsilon(W)-V_0(W)\right]\right]^2.  
\end{multline*}
We deduce from proposition \ref{th:erange} that for every
$\varepsilon\in (0,\varepsilon_0]$,
\[
\E^X_\nu\left[\left[\sum_{i\in I}
V_\varepsilon(W^i)-\sum_{i\in I}
V_0(W^i)\right]^2\right]
\leq  [(\nu,\ind)+(\nu,\ind)^2]h_d(\varepsilon)^\kappa.
\]

\noindent
Notice the hypothesis on
$\supp \nu$  and $\theta$ imply that $\lim_{\varepsilon\rightarrow 0} \varphi_d(\varepsilon)\val{
  (\supp \nu )^{h_d(\varepsilon)}}=0$.
Arguments similar to those used in the first part of the proof of corollary
\ref{cor:cvpsrangeY}  show then  a.s.
\[
\lim_{\varepsilon\rightarrow 0} \sum_{i\in I}
V_\varepsilon(W^i)=\sum_{i\in I}
V_0(W^i).
\]
Notice we have $\sum_{i\in I}
V_0(W^i)=\cstap \int_0^\infty  ds \; (X_s,\ind_A)$. Using the above remark on $\supp \nu$, we deduce that a.s.
\[
\limsup_{\varepsilon\rightarrow 0} \varphi_d(\varepsilon)
\val{\range_0(X)^{\varepsilon} \cap A} \leq  \cstap\int_0^\infty
  ds\;(X_s,\ind_A).
\]

Second step. To get a lower bound, consider an increasing sequence $(E_p,p\geq 1)$ of
measurable subsets of $E=C(\R^+,\cw)$ such that $\bigcup_{p\geq 1}
E_p=E$ and 
$\int \nu(dx) \N_x [E_p]=\alpha_p<\infty $. (For instance we can take
$E_p=\left\{W; \sup_{s\geq 0} \zeta_s\geq 1/p\right\}$.) Then a.s. the set
$I_p=\left\{i\in I; W^i \in E_p\right\}$ is finite. We have 
\[
\varphi_d(\varepsilon)\val{\range_0(X)^{\varepsilon} \cap A}\geq \sum_{i\in
  I_p} V_\varepsilon(W^i) -\sum_{(i,j)\in I_p^2;\; i \neq j}
  U_\varepsilon(W^i,W^j),
\]
where
\begin{align*}
U_\varepsilon(W^i,W^j)
&=\varphi_d(\varepsilon) \val{\range(W^i)^{\varepsilon}
  \cap \range(W^j)^{\varepsilon}\cap A \cap \bar B(x^i_0,
  h_d(\varepsilon))^c \cap \bar B(x^j_0,
  h_d(\varepsilon))^c}\\
&=\varphi_d(\varepsilon)\int_{A \cap \bar B(x^i_0,
  h_d(\varepsilon))^c \cap \bar B(x^j_0,
  h_d(\varepsilon))^c}dy\;
  \ind_{\left\{T_{(y,\varepsilon)}(W^i)<\infty
  \right\}}\ind_{\left\{T_{(y,\varepsilon)}(W^j)<\infty \right\}} .
\end{align*}
Arguments similar to those of the first step show that a.s.
\[
\lim_{\varepsilon\rightarrow 0} \sum_{i\in I_p} V_\varepsilon(W^i)
  =\sum_{i\in I_p} V_0(W^i)= \sum_{i\in I_p}  \cstap\int_0^\infty 
  ds\;(Y_s(W^i),\ind_A).
\]
Now conditionally on the cardinal of $I_p$, the Brownian snakes $(W^i, i\in
I_p)$ are independent and have the same law: $\mu_p=\alpha_p^{-1}\int \nu(dx) \N_x
[\cdot\cap E_p]$. For two independent Brownian snakes
$(W,W')$ under $\mu_p\otimes \mu_p$, we get using  
\reff{eq:umajo}, that for $\varepsilon$ satisfying (C),
\begin{align*}
   \mu_p\otimes \mu_p[U_\varepsilon(W,W')]
   & \leq  \alpha_p^{-2} \iint \nu(dx_0)  \nu(dx'_0)
   \N_{x_0}\otimes \N_{x'_0} [U_\varepsilon(W,W')]\\
   & \leq  \varphi_d(\varepsilon)\alpha_p^{-2} \iint \nu(dx_0)  \nu(dx'_0)
   \int_{A \cap \bar B(x_0,
   h_d(\varepsilon))^c \cap \bar B(x'_0,
   h_d(\varepsilon))^c}dy\;\\
   &\hspace{4cm}\left[\cstb\varphi_d(\varepsilon)^{-1}\val{y-x_0}^{2-d}\right]
   \left[\cstb\varphi_d(\varepsilon)^{-1}\val{y-{x'_0}}^{2-d}\right]\\
   & \leq  \varphi_d(\varepsilon)^{-1}\alpha_p^{-2} (\nu,\ind)^2 \cstb^2
   \sup_{x_0\in \R^d} 
   \int_{\bar B(0,R_0) \backslash \bar  B(x_0, h_d(\varepsilon))} dy
   \val{y-x_0}^{4-2d}\\
   & \leq  \cst \varphi_d(\varepsilon)^{-1} h_d(\varepsilon)^{4-d} \quad \text{if $d\geq 5$}\\
   & \leq  \cst \varphi_4(\varepsilon)^{-1} \log(\log(1/\varepsilon))\quad \text{if $d=4$}\\
   & \leq  \cst h_d(\varepsilon)^{\theta/2} \quad \text{if $d\geq 4$},
\end{align*}
where the constant $\cst$ is independent of $\varepsilon$ and $A$. Using the
Borel-Cantelli lemma for the sequence
$(h_d(\varepsilon_n)=n^{-4/\theta},n\geq 1)$, and a monotonicity
argument, we get that $\mu_p\otimes \mu_p$-a.s. 
$\lim_{\varepsilon\rightarrow 0} U_{\varepsilon}(W,W')=0$. Then since
the cardinal of $I_p$ is 
a.s. finite, we get that for every integer $p\geq 1$, a.s.,
\[
\lim_{\varepsilon\rightarrow 0}
\sum_{(i,j)\in I_p^2;\; i \neq j}
  U_\varepsilon(W^i,W^j)=0.
\]
We deduce that for every integer $p\geq 1$, a.s.
\[
\liminf_{\varepsilon\rightarrow 0} \varphi_d(\varepsilon)
\val{\range_0(X)^\varepsilon\cap A} \geq \sum_{i\in I_p}
\cstap\int_0^\infty  ds\; (Y_s(W^i), \ind_A).
\]
We get the lower bound by letting $p\rightarrow \infty $. This and the upper bound of the first step ends the proof of the theorem.
\findemo

\section{Proof of proposition \ref{th:erange}}
\label{sec:prprop}
We shall use many times in the sequel the fact that $\int_0^\infty  ds\; (Y_s,\ind_A)=\int_0^\sigma ds \;
\ind_A(\hat W_s)$ $\N_{x_0}$-a.e.  We assume $d\geq 4$. 
We recall easy equalities, which can  readily be deduced from the
results of section \ref{sec:annexemo}. For every $A\in
\cb(\R^d)$,  we have
\begin{align}
   \label{eq:N11}
   \N_x\left[\int_0^\sigma ds\;\ind_A(\hat W_s)
   \right]
   &=\int_A dy\; G(x,y),\\
   \intertext{where $G$ is the Green kernel in $\R^d$:
     $G(x,y)=2^{-1}\pi^{-d/2}\Gamma([d-2]/2)\val{x-y}^{2-d}$, and}
   \label{eq:N12}
   \N_x\left[\left[\int_0^\sigma ds\;\ind_A(\hat W_s)
   \right]^2\right]
   &= 4\int dy\; G(x,y)\left[\int_A dz\; G(y,z)\right]^2.
\end{align}
We can also compute the first moment under $\ce^*_\rw$. For every $A\in
\cb(\R^d)$, $\rw \in \cw$, we have with $\zeta=\zeta_{(\rw)}$,
\begin{equation}
\label{eq:Egreen}
   \hspace{-1.5cm}
   \ce^*_\rw
   \left[\int_0^\sigma ds\;\ind_A(\hat W_s)\right]= 2\int_0^\zeta dt\;
   \N_{\rw(t)}\left[\int_0^\sigma ds\;\ind_A(\hat W_s) \right]=
   2\int_0^\zeta dt \int_A dy\; G(\rw(t),y).
\end{equation}

\noindent
Thanks to the space invariance
of the law of the Brownian snake, we shall only consider  the case
$x_0=0$ and $A\subset \bar B(0,R_0)$, for $R_0$ fixed.  
We fix $\theta\in (0,1/d)$ and $R_0>1$. Let
$\varepsilon'_0>0$ satisfying (C). We  consider $\varepsilon\in
(0,\varepsilon'_0)$. In this section, we  denote by $\cst$,
$c_1,c_2,\ldots$ positive constants whose values depend only on $d,\theta$ and
$R_0$. The value of $\cst$ may vary from line to line.
For short we shall write
$A_\varepsilon=A\cap \bar B(0,h_d(\varepsilon))^c  $ (not to be
confused with $A^\varepsilon$) and $\range$ for $\range(W)$. 

We first consider the case $d\geq 5$.
Notice that 
\[
\N_0\left[\val{\range^\varepsilon\cap A_\varepsilon}\right]
=\int_{A_\varepsilon} dx \;\N_0\left[T_{(x,\varepsilon)}<\infty \right]
=\int_{A_\varepsilon} dx\;u_\varepsilon(x).
\]
Thus we deduce from \reff{eq:umino} and
\reff{eq:udel}, that for $\varepsilon\in (0,\varepsilon'_0)$,
\begin{multline*}
\csta \varepsilon^{d-4} \int_A dx\; \val{x}^{2-d} - \csta
\varepsilon^{d-4} \int_{B(0, \varepsilon^{1-\theta})} dx \;
\val{x}^{2-d}\\
%\begin{aligned}
 \leq \N_0\left[\val{\range^\varepsilon\cap A_\varepsilon}\right]
\leq  \varepsilon^{d-4}[\csta +\cstaa h_d(\varepsilon)^{\theta/2}] \int_A dx\; \val{x}^{2-d}.     
%\end{aligned}
\end{multline*}
Therefore using also \reff{eq:N11}, we have 
\[
  \val{\N_{x_0}\left[\varepsilon^{4-d}\val{\range(W)^\varepsilon\cap A_\varepsilon}-\cstap \int_0^\infty
  ds\;(Y_s,\ind_A)\right]}
\leq \cst h_d(\varepsilon)^{\theta/2}.
\]
Thus we get the first bound  of proposition \ref{th:erange} (take
$\kappa<\theta$ and $\varepsilon_0$ small enough). The proof is similar for $d=4$ (use \reff{eq:udel41} instead of \reff{eq:udel1} and the fact that $\val{x}$ is bounded by $R_0$).

Now we  will prove the second bound. To this end
we have to find an upper bound on $I=\N_0\left[\val{\range^\varepsilon\cap A_\varepsilon}^2\right]$ and a
lower bound on $J=\N_0\left[\val{\range^\varepsilon\cap
A_\varepsilon}\int_0^\sigma ds\; \ind_A(\hat W_s) \right]$.

\subsection{An upper bound on $I$}
\label{sec:majoI}

The term $I$ can also be written 
\[
I=\iint_{A_\varepsilon\times A_\varepsilon} dx\;dy \; \N_0 \left[
T_{(x,\varepsilon)}<\infty ; T_{(y,\varepsilon)}<\infty \right].
\] 
Consider the above integral as the sum of the integral over 
$\val{x-y}\leq 2h_d(\varepsilon)$ (denoted by $I_1$) and the one
over 
$\val{x-y}> 2h_d(\varepsilon)$ (denoted by $I_2$). Using \reff{eq:umajo} we  easily obtain an upper bound on $I_1$: 
\initre
\begin{align*}\initer
      I_1
      & \leq \val{B(0,2h_d(\varepsilon))} \int_{A_\varepsilon} dx\; \N_0 \left[
        T_{(x,\varepsilon)}<\infty \right]\\
      & \leq \cst h_d(\varepsilon)^d \int_{A_\varepsilon} dx\;
      \varphi_d(\varepsilon)^{-1}
      \cstb \val{x}^{2-d} \leq \cste   \varphi_d(\varepsilon)^{-2}h_d(\varepsilon)^3.
\end{align*}
\setcounter{cstei}{\thecste}

\noindent
Notice  the event % the set
$\left\{T_{(x,\varepsilon)}<\infty ; T_{(y,\varepsilon)}<\infty
\right\}$ is a subset of 
\[
\left\{T_{(x,\varepsilon)}<\infty ; T_{(y,\varepsilon)}\circ
\theta_{T_{(x,\varepsilon)}}<\infty \right\}\cup
\left\{T_{(y,\varepsilon)}<\infty ; T_{(x,\varepsilon)}\circ 
\theta_{T_{(y,\varepsilon)}}<\infty \right\},
\]
where $\theta_t$ is the usual 
shift operator. By symmetry, we get
\begin{equation}
   \label{eq:I2majo}
\hspace{-2cm}I_2\leq 2\iint_{A_\varepsilon\times A_\varepsilon}dx\;dy\;
\ind_{\left\{\val{x-y}>2h_d(\varepsilon)\right\}} \N_0
\left[T_{(x,\varepsilon)}<\infty ; T_{(y,\varepsilon)}\circ 
\theta_{T_{(x,\varepsilon)}}<\infty\right].
\end{equation}
Using the strong Markov property of the Brownian snake under $\N_0$ at
the 
stopping time $T_{(x,\varepsilon)}$ and (\ref{eq:hitting}), we see that
the quantity $\N_0\left[T_{(x,\varepsilon)}<\infty ; T_{(y,\varepsilon)}\circ 
\theta_{T_{(x,\varepsilon)}}<\infty\right]$  is equal to
\[
\N_0
\left[T_{(x,\varepsilon)}<\infty ; 
2\int_0^{\zeta_{T_{(x,\varepsilon)}}\wedge
  \tau_{B(y,\varepsilon)}(W_{T(x,\varepsilon)})} 
dt\;u_\varepsilon\left(W_{T_{(x,\varepsilon)}}(t)-y\right)  
\expp{\left[-2\int_0^t u_\varepsilon\left(W_{T_{(x,\varepsilon)}}(s)-y\right) ds
  \right]}\right] .
\]
Finally the law of the stopped path $W_{T_{(x,\varepsilon)}}$ under
$\N_0$ is given by (\ref{eq:Nbeta}). Thus the previous  expression is equal
to
\[
2\int_0^\infty  \!dt\;\E_0
\left[\tau_{B(x,\varepsilon)}>t;\tau_{B(y,\varepsilon)}>t; 
u_\varepsilon(\beta_t-x)u_\varepsilon(\beta_t-y)
\expp{\left[-2\int_0^t ds\;
\left[u_\varepsilon(\beta_s-x)+u_\varepsilon(\beta_s-y)
\right]\right]} \right].
\]
We substitute this last expression for
$\N_0\left[T_{(x,\varepsilon)}<\infty ; T_{(y,\varepsilon)}\circ  
  \theta_{T_{(x,\varepsilon)}}<\infty\right]$  in (\ref{eq:I2majo}), and
then
decompose the right-hand side of (\ref{eq:I2majo}) in three terms by
considering the integral in $dxdy$ on the sets 
$\val{\beta_t-x}\wedge\val{\beta_t-y}>h_d(\varepsilon)$ (integral
$I_{21}$), 
$\val{\beta_t-x}\leq h_d(\varepsilon)$ (integral
$I_{22}$), and 
$\val{\beta_t-y}\leq h_d(\varepsilon)$ (integral
$I_{23}$) (recall $\val{x-y}>2h_d(\varepsilon)$).

\textbf{An upper bound on $I_{21}$.}
We shall need the following notation:
\[
I_0=4 \csta^2 
      \int dz \; G(0,z) \left[\int_{A}dx\;
      \val{z-x}^{2-d}\right]^2.
\]
We use \reff{eq:udel} to bound $I_{21}$ above by: for $\varepsilon\in (0,\varepsilon'_0)$,
\initre
\begin{align*}\initer
      &4\iint_{A_\varepsilon\times A_\varepsilon}dx\;dy\;
      \ind_{\{\val{x-y}>2h_d(\varepsilon)\}} 
      \int_0^\infty dt\;
      \E_0\Big[ \val{\beta_t-x}>h_d(\varepsilon);
      \val{\beta_t-y}>h_d(\varepsilon) ;\\
      &\hspace{1cm}\varphi_d(\varepsilon)^{-2}\val{ \beta_t-x }^{2-d}
      \val{\beta_t-y}^{2-d} \left(\csta +\cstaa h_d(\varepsilon)^{\theta/2}\right)^2 \Big]
      \\
      &\hspace{.5cm} \leq 4 \varphi_d(\varepsilon)^{-2} \left[\csta^2 +\cst h_d(\varepsilon)^{\theta/2}\right]
      \iint_{A\times A}dx\;dy\; \int dz \; G(0,z)
      \val{z-x}^{2-d}\val{z-y}^{2-d}\\
      &\hspace{.5cm} \leq \varphi_d(\varepsilon)^{-2} I_0+\cste \varphi_d(\varepsilon)^{-2}h_d(\varepsilon)^{\theta/2}.
\end{align*}
\setcounter{csteii}{\thecste}

\textbf{An upper bound on $I_{22}$ and $I_{23}$.}
By symmetry we have $I_{22}=I_{23}$. Before getting an upper bound on
$I_{22}$, notice that $\val{\beta_t-x}\leq h_d(\varepsilon)$ and
$\val{x-y}>2h_d(\varepsilon)$  imply
$\val{\beta_t-y}>h_d(\varepsilon)$. Furthermore thanks to
\reff{eq:umajo}, we get
\initre
\begin{align*}\initer
      \int_{A_\varepsilon}
      dy\;\ind_{\left\{\val{\beta_t-y}>h_d(\varepsilon)\right\}}
      u_\varepsilon(\beta_t-y) \expp{-2\int_0^t u_\varepsilon(\beta_s-y)
      ds}
      & \leq \int_{A}
      dy\;\left[\cstb \varphi_d(\varepsilon)^{-1}\val{\beta_t-y}^{2-d}
      \right]
      \\
      &\leq \cstb \varphi_d(\varepsilon)^{-1} \int_{B(0,R_0)} dy\; \val{y}^{2-d} =\cste \varphi_d(\varepsilon)^{-1}.
\end{align*}
Thus the sum $I_{22}+I_{23}$ is bounded  above by
\[
8\mcste  \varphi_d(\varepsilon)^{-1}\int_{A_\varepsilon}dx\; \int_0^\infty
dt\; \E_0 \left[\tau_{B(x,\varepsilon)}>t;
\ind_{\left\{\val{\beta_t-x}\leq h_d(\varepsilon)\right\}}
u_\varepsilon(\beta_t-x) \expp{-2\int_0^t u_\varepsilon(\beta_s-x)ds } 
\right].
\]
Using the Cauchy-Schwarz inequality and formula \reff{eq:fk}, we
get
\initre
\begin{align*}\initer
      I_{22}+I_{23}
      &\leq  8\mcste  \varphi_d(\varepsilon)^{-1}
      \left[\int_{A_\varepsilon}dx\; \int_0^\infty dt \;\P_0 \left[
      {\val{\beta_t-x}\leq h_d(\varepsilon)}\right]\right]^{1/2}
      \\
      &\phantom{\leq  2\mcste R_0^2 \varepsilon^{d-4} }
      \times \left[\int_{A_\varepsilon}dx\; \int_0^\infty dt\; \E_0 \left[
      \tau_{B(x,\varepsilon)}>t;
      u_\varepsilon(\beta_t-x)^2 \expp{-4\int_0^t
      u_\varepsilon(\beta_s-x)ds }\right]\right]^{1/2} 
      \\
      &\leq  8\mcste  \varphi_d(\varepsilon)^{-1}
      \left[\int_{A}dx\; \int dz \;G(0,z) \ind_{
      \left\{\val{z-x}\leq h_d(\varepsilon)\right\}}\right]^{1/2}
      \left[\int_{A_\varepsilon}dx\; 2^{-1} u_\varepsilon(x) \right]^{1/2}.
\end{align*}
Then thanks to \reff{eq:umajo},
we get $I_{22}+I_{23}
\leq \cste \varphi_d(\varepsilon)^{-3/2}h_d(\varepsilon)^{d/2}\leq \mcste \varphi_d(\varepsilon)^{-2}h_d(\varepsilon)^{3/2} $. 
\setcounter{csteiii}{\thecste}

\textbf{Conclusion on the upper bound on $I$.}
By combining the previous results, we get for $d\geq 4$
\[
I\leq \cstei  \varphi_d(\varepsilon)^{-2}h_d(\varepsilon)^{3}
+\varphi_d(\varepsilon)^{-2} I_0
+\csteii \varphi_d(\varepsilon)^{-2}h_d(\varepsilon)^{\theta/2}
+\csteiii \varphi_d(\varepsilon)^{-2}h_d(\varepsilon)^{3/2}. 
\]
Thus  we get
$ \varphi_d(\varepsilon)^2 I\leq I_0+\cste h_d(\varepsilon)^{\theta/2}$.
\setcounter{csteiii}{\thecste}

\subsection{A lower bound on $J$}
\label{sec:sublbJ}

We shall need the last hitting time of $\bar B(x,\varepsilon)$ under $\N_0$
for the Brownian snake:
\[
L_{(x,\varepsilon)}=\sup \left\{s\geq 0; \exists t\in [0,\zeta_s],
  W_s(t)\in \bar B(x,\varepsilon)\right\}.
\]
We then get
\initre
\begin{multline*}\initer
J=\int_{A_\varepsilon} dx \;\N_0 \left[T_{(x,\varepsilon)}<\infty ;
  \int_0^{L_{(x,\varepsilon)}} ds \;\ind_A(\hat W_s) \right]\\
+\int_{A_\varepsilon} dx \;\N_0 \left[T_{(x,\varepsilon)}<\infty ;
  \int_{T_{(x,\varepsilon)}}^\sigma ds \;\ind_A(\hat W_s) \right]\\
-\int_{A_\varepsilon} dx \;\N_0 \left[T_{(x,\varepsilon)}<\infty ;
  \int_{T_{(x,\varepsilon)}}^{L_{(x,\varepsilon)}} ds \;\ind_A(\hat W_s)
  \right] .
\end{multline*}
The time-reversal invariance property of the It\^o measure and the
characterization of the excursion measure $\N_x$ readily imply that the
latter itself enjoys the same invariance property. Thus the 
first two terms of the right-hand side are equal. We shall denote their
sum by $J_1$. Let $J_2$ denote the third term.

\textbf{A lower bound on $J_1$.}
Let us use the strong Markov property of the Brownian snake at time
$T_{(x,\varepsilon)}$,  then (\ref{eq:Egreen}) and (\ref{eq:Nbeta}), to
get 
\initre
\begin{align*}\initer
      J_1
      &= 2 \int_{A_\varepsilon} dx \;\N_0 \left[T_{(x,\varepsilon)}<\infty ;
        2 \int_0^{\zeta_{T_{(x,\varepsilon)}}} dt \; \int_A dy \;
      G\left(W_{T_{(x,\varepsilon)}} (t),y\right) \right]\\
      &= 4 \int_{A_\varepsilon} dx \int_A dy \int_0^\infty  dt \;\E_0
      \left[\tau_{B(x,\varepsilon)}> t ;
      G(\beta_t,y) u_\varepsilon(\beta_t-x)\expp{-2\int_0^t
      u_\varepsilon(\beta_s-x)ds}\right].
\end{align*}
Fatou's lemma gives that $ \liminf_{\varepsilon\rightarrow 0}
\varphi_d(\varepsilon) J_1\geq  J_0$, where 
\[
J_0=4\csta \iint_{A\times A} dxdy\int dz \; G(0,z)G(z,y)
\val{z-x}^{2-d}.
\]
Unfortunately, we need an estimate on the rate of convergence. This
requires  some technical calculations. Notice  that on 
$\{\tau_{B(x,h_d(\varepsilon))}(\beta)>t\}$, inequalities
  \reff{eq:umino}, \reff{eq:umino4} and \reff{eq:umajo} imply
\[
     \csta \varphi_d(\varepsilon)^{-1} F_d(\beta_t-x)
    \val{\beta_t-x}^{2-d}\leq u_\varepsilon(\beta_t -x)\leq\cstb \varphi_d(\varepsilon)^{-1}
    \val{\beta_t-x}^{2-d},
\]
where $F_d(z)=1$ if $d\geq 5$ and $F_4(z)=\left[1+\log (2\val{z})/\log(1/\varepsilon)\right]^{-1}$.
For short we write $\Gamma_t=2\cstb \varphi_d(\varepsilon)^{-1} \int_0^t
\val{\beta_s-x}^{2-d} ds$. Then $\varphi_d(\varepsilon)J_1$ is bounded below by
\[
J'_1=4\csta  \int_{A_\varepsilon}\!dx \int_A \!dy
\int_0^\infty \!dt \;\E_0\left[\tau_{B(x,h_d(\varepsilon))}>t; G(\beta_t,y)
\val{\beta_t-x}^{2-d}F_d(\beta_t-x)\expp{-\Gamma_t }\right].
\]
In order to obtain an upper bound on $\val{ J'_1-J_0}$,
we have to find an upper bound on
\[
 \iint_{A\times A}dx  dy
\int_0^\infty dt \;\E_0\left[ G(\beta_t,y)\val{\beta_t-x}^{2-d}
  \left[1-\ind_{A_\varepsilon}(x) \ind_{ \left\{\tau_{B(x,h_d(\varepsilon))}>t\right\}} 
F_d(\beta_t-x)\expp{-\Gamma_t }\right]\right].
\]
Thus we shall decompose $1-\ind_{A_\varepsilon}(x) \ind_{\left\{\displaystyle
\tau_{B(x,h_d(\varepsilon))}>t\right\}} F_d(\beta_t-x)
\expp{-\Gamma_t }$ into a sum of four  terms:
\begin{multline*}
      \left[1-\ind_{A_\varepsilon}(x)\right]
      +\ind_{A_\varepsilon}(x)\left[1-\ind_{\left\{
      \tau_{B(x,h_d(\varepsilon))}>t\right\}} \right] \\
      +\ind_{A_\varepsilon}(x)\ind_{\left\{
      \tau_{B(x,h_d(\varepsilon))}>t\right\}}\left[1-F_d(\beta_t-x) \right]
      +\ind_{A_\varepsilon}(x)\ind_{\left\{
      \tau_{B(x,h_d(\varepsilon))}>t\right\}}F_d(\beta_t-x)\left[1-\expp{-\Gamma_t } \right]  .
\end{multline*}
We denote by  $J_{11}$, $J_{12}$, $J_{13}$ and $J_{14}$ the corresponding
integrals. The integral 
\[
J_{11}= \int_{A\backslash A_\varepsilon}dx
\int_A  dy
\int_0^\infty dt \;\E_0\left[ G(\beta_t,y)\val{\beta_t-x}^{2-d} \right]
\]
is easily bounded above by
\[
  \int_{B(0,h_d(\varepsilon))} dx\; \int_{B(0,R_0)} dy \int dz \;
G(0,z) G(z,y)\val{z-x}^{2-d}
\leq \cste h_d(\varepsilon)^{2}.
\]
We bound  $J_{12}$ by applying the strong Markov property of  Brownian motion at time
$\tau_{B(x,h_d(\varepsilon))}$,
\initre\setcounter{csteij}{\thecste}
\begin{align*}\initer
     J_{12}&=\int_{A_\varepsilon}dx \int_A dy
     \int_0^\infty dt \;\E_0\left[\tau_{B(x,h_d(\varepsilon))}\leq t ; G(\beta_t,y)
     \val{\beta_t-x}^{2-d}\right]\\
     %&\hspace{1cm} 
     &\leq \int_{A_\varepsilon}dx \int_A dy\;
     \E_0\left[\tau_{B(x,h_d(\varepsilon))}<\infty ; \int
     dz \;G(\beta_{\tau_{B(x,h_d(\varepsilon))}},z) G(z,y)
     \val{z-x}^{2-d}\right].
\end{align*}    
An easy calculation shows that there exists a constant $\cste$ such that
for every  $(x,x')\in B(0,2R_0)\times B(0,2R_0)$, $\val{x-x'}\leq 1/2$,
\[
\int_{B(0,R_0)}dy \int dz \;G(x',z)G(z,y) \val{z-x}^{2-d}
\leq 
\mcste \varphi_d(\val{x'-x})
\]
Furthermore %thanks to proposition 1.6 p56 from \cite{ps:bm}, we get for
we have for every $r\in (0,1)$, 
\begin{equation}
   \label{eq:intta}
   \int_{B(0,R_0)}dx \;\P_0\left[\tau_{B(x,r)}<\infty  \right]
   =\int_{B(0,R_0)}dx\; \left[\left(\frac{r}{\val{x}}\right)^{d-2}\wedge 1\right]
   \leq \cst r^{d-2}.
\end{equation}
We deduce from the
previous remarks that if $d\geq 5$,
\initre
\[
\initer
   J_{12}
   \leq \cst  h_d(\varepsilon)^{4-d}
   \int_{A_\varepsilon} dx\; \P_0\left[
   \tau_{B(x,h_d(\varepsilon))}<\infty \right]
   \leq \cst h_d(\varepsilon)^{4-d+d-2}=\cst h_d(\varepsilon)^2, 
\]
and if $d=4$, $J_{12}\leq \cst \log(1/h_d(\varepsilon)) h_d(\varepsilon)^2 $. Thus we get that for $d\geq 4$, $J_{12}\leq \cste h_d(\varepsilon)^{3/2}$.
\setcounter{cstej}{\thecste}

\noindent
If $d\geq 5$ then $J_{13}=0$. For $d=4$ thanks to \reff{eq:lolo} we have for $\val{z}\geq h_4(\varepsilon)$, $\val{1-F_4(z)}\leq 2 \val{\log(2\val{z})}/\log(1/\varepsilon)$. We deduce that 
\initre
\begin{align*}\initer
   J_{13}
   &\leq \log(1/\varepsilon)^{-1} \iint_{A\times A}\!\!dx  dy \int_0^\infty \!dt\; \E_0\left[\tau_{B(x,h_d(\varepsilon))}>t; G(\beta_t,y)2\val{\log({2\val{\beta_t-x}})} \val{\beta_t-x}^{-2}\right]\\
   &\leq \cst \log(1/\varepsilon)^{-1} \iint_{A\times A}dx  dy \int dz\; G(0,z)\val{\log({2\val{z-x}})} \val{z-x}^{-2}G(z,y)\\
   &\leq \cst \log(1/\varepsilon)^{-1} \leq \cste h_d(\varepsilon)^\theta.
\end{align*}

\noindent
Notice first that thanks to \reff{eq:lolo}, $F_d(z)\leq 2$ for $\val{z}\geq h_d(\varepsilon)$. We have, using the Markov property for Brownian motion at time
$s$, 
\setcounter{cstejjiii}{\thecste}
\initre
\begin{align*}\initer
   J_{14}
   &\leq 2\iint_{A\times A}dx  dy
    \int_0^\infty dt\; \\
&\hspace{2cm}\E_0\Bigg[
     \tau_{B(x,h_d(\varepsilon))}>t;
    G(\beta_t,y)
     \val{\beta_t-x}^{2-d}2\cstb \varphi_d(\varepsilon)^{-1}\int_0^t
     \val{\beta_s-x}^{2-d} ds\Bigg]\\ 
   &\leq \cst\varphi_d(\varepsilon)^{-1}\iint_{A\times A} dxdy\;\int_0^\infty  ds\; \int_0^\infty  dt\;\\
&\hspace{2cm}
   \E_{0}\left[\val{\beta_s-x}^{2-d}\E_{\beta_s}
   \left[\val{\beta_{t}-x}>h_d(\varepsilon) ;
   G(\beta_t,y)\val{\beta_{t}-x}^{2-d}\right]\right]\\
   &\leq \cst \varphi_d(\varepsilon)^{-1} M(d,h_d(\varepsilon)),
\end{align*}
where
\[
M(d,\varepsilon)=\iint_{B(0,R_0)^2} dxdy\;\iint dzdz'\;
   G(0,z)\val{z-x}^{2-d} G(z,z') G(z',y)\val{z'-x}^{2-d}\ind_{\val{z'-x}>\varepsilon}.
\]
An easy computation shows there exists a constant $c$ such that for $\varepsilon\in (0,1]$, 
\begin{equation}
   \label{eq:Mde}
M(d,\varepsilon)\leq \left\{ 
\begin{array}{ll}
   \cst & \text{ if } d\in \{4,5\},\\
   \cst+ \cst \log(1/\varepsilon) & \text{ if } d=6,\\
   \cst \varepsilon^{6-d} & \text{ if } d\geq 7.
\end{array}
\right.
\end{equation}
Thus we easily deduce that $J_{14}\leq \cste h_d(\varepsilon)^\theta$.

We have  $\varepsilon^{4-d}J_1\geq J_0-4\csta
(J_{11}+J_{12}+J_{13}+J_{14})$.
Putting  together the previous results, we get for $d\geq 4$,
\setcounter{csteji}{\thecste}
\[
   \varepsilon^{4-d}J_1\geq 
    J_0- 4\csta [\csteij h_d(\varepsilon)^2+\cstej  h_d(\varepsilon)^{3/2}
+\cstejjiii h_d(\varepsilon)^\theta   
+\csteji 
   h_d(\varepsilon)^\theta]\geq 
J_0- \cste h_d(\varepsilon)^\theta.
\]
\setcounter{cstejii}{\thecste}

\textbf{An upper bound on $J_2$.}
We will first recall the decomposition of the Brownian snake under $\ce^*_{\rw}$ (see theorem 2.5 in \cite{lg:pvmppde}). We denote by $(\alpha_i,\beta_i)$, $i\in I$, the excursion intervals of $\zeta$ above its minimum process (i.e. of the process $(\zeta_t-\inf_{s\in [0,t]}\zeta_s)$ above $0$) before $\sigma$ under $\ce^*_{\rw}$. For $i\in I$ the paths $W_s, s\in [\alpha_i,\beta_i]$ coincide over $[0,\zeta_{\alpha_i}]$. For every $i\in I$, and $s\geq 0$ we set $W_s^i(t)=W_{(\alpha_i+s)\wedge \beta_i}(t+\zeta_{\alpha_i})$, $t\in [0,\zeta^i_s]$ with $\zeta^i_s=\zeta_{(\alpha_i+s)\wedge \beta_i}-\zeta_{\alpha_i}$. Then $W^i_s$ is a stopped path ($W^i_s\in \cw$) with initial point $W_{(\alpha_i+s)\wedge \beta_i}(\zeta_{\alpha_i})=\hat W_{\alpha_i}=\rw(\zeta_{\alpha_i})$. 
\begin{proposition}[Le Gall]
The random measure $\sum_{i\in I} \delta_{(\zeta_{\alpha_i},W^i)}$ is under $\ce^*_{\rw}$ a Poisson point measure on $[0,\zeta_{(\rw)}]\times C(\R^+,\cw)$ with intensity $2 dt\; \N_{\rw(t)}[\cdot]$.
\end{proposition}
The process $\sum_{i\in I, \zeta_{\alpha_i}\leq t} \delta_{W^i}$ for $t\in [0,\zeta_{\rw}]$ is a Poisson point process with inhomogeneous intensity. We will now describe the law under $\ce^*_{W_{T_{(x,\varepsilon)}}}$ of the first excursion $(\zeta_{\alpha_{i_0}}, W^{i_0})$ which hits the ball $\bar B(x,\varepsilon)$, that is, with evident notation, the first excursion for which $T_{(x,\varepsilon)}(W^i)$ is finite. We first notice that under $\N_0[\cdot\mid T_{(x,\varepsilon)}<\infty ]$, $\ce^*_{W_{T_{(x,\varepsilon)}}}$-a.s. there are such excursions. Indeed we have thanks to lemma 2.1 of \cite{lg:hppt} that $\N_0[.\mid T_{(x,\varepsilon)}<\infty ]$-a.s.
%\begin{align*}
\[
  \ce^*_{W_{T_{(x,\varepsilon)}}}[\exists i\in I, T_{(x,\varepsilon)}(W^i)<\infty ]
%&=1-\exp{-2 \int_0^{\zeta_{T_{(x,\varepsilon)}}} dt\; \N_{W_{T_{(x,\varepsilon)}}(t)}[T_{(x,\varepsilon)}<\infty ]}\\
=1-\exp{-2 \int_0^{\zeta_{T_{(x,\varepsilon)}}} dt\; u_\varepsilon(W_{T_{(x,\varepsilon)}}(t)-x)}=1.
\]
%\end{align*}
Since the integral $\int_0^r  dt\; u_\varepsilon(W_{T_{(x,\varepsilon)}}(t)-x)$ is  finite for $r<\zeta_{T_{(x,\varepsilon)}}$, we deduce there exists a unique first excursion $i_0$ which hits $\bar B(x,\varepsilon)$. Classical arguments on Poisson point process implies that the law of $(\zeta_{\alpha_{i_0}}, W^{i_0})$ is $2 \ind_{[0,\zeta_{{T_{(x,\varepsilon)}}})}(t)dt\; \N_{W_{T_{(x,\varepsilon)}}(t)}[T_{(x,\varepsilon)}<\infty ,\cdot]$. We introduce the random time $M_{(x,\varepsilon)}=\inf \left\{s>T_{(x,\varepsilon)}; \zeta_s=m(T_{(x,\varepsilon)}, L_{(x,\varepsilon)})\right\}$. It is clear from the definition of the excursion $i_0$ that ${\alpha_{i_0}}=M_{(x,\varepsilon)}$ under $\ce^*_{W_{T_{(x,\varepsilon)}}}$. We will now express $J_2$ using the excursion $i_0$.
We have
\begin{align*}
   J_2
&=2\int_{A_\varepsilon} dx\;  \N_0\left[T_{(x,\varepsilon)}<\infty ; \int_{M_{(x,\varepsilon)}}^{L_{(x,\varepsilon)}}ds\;\ind_A (\hat W_s)\right]\\
&=2\int_{A_\varepsilon} dx\;  \N_0\left[T_{(x,\varepsilon)}<\infty ; \ce^*_{W_{T_{(x,\varepsilon)}}}\left[\int_{M_{(x,\varepsilon)}}^{L_{(x,\varepsilon)}}ds\;\ind_A (\hat W_s)\right]\right]\\
&=2\int_{A_\varepsilon} dx\;  \N_0\left[T_{(x,\varepsilon)}<\infty ; \ce^*_{W_{T_{(x,\varepsilon)}}}\left[\int_{\alpha_{i_0}}^{L_{(x,\varepsilon)}(W^{i_0})}ds\;\ind_A (\hat W_s^{i_0})\right]\right]\\
&=4\int_{A_\varepsilon} dx\;  \N_0\left[T_{(x,\varepsilon)}<\infty ; \int_0^{\zeta_{{T_{(x,\varepsilon)}}}}dt\; \N_{W_{T_{(x,\varepsilon)}}(t)}\left[T_{(x,\varepsilon)}<\infty ; \int_{0}^{L_{(x,\varepsilon)}}ds\;\ind_A (\hat W_s)\right]\right].
\end{align*}
We used the time reversal property of the Brownian snake for the first equality, then the strong Markov property and at last the definition of the excursion $i_0$ and its law. We will distinguish according to $\left\{t\geq \tau_{B(x,h_d(\varepsilon))}\right\}$ (integral $J_{21}$) and $\left\{t<\tau_{B(x,h_d(\varepsilon))}\right\}$
 (integral $J_{22}$). Notice that since $x\in A_\varepsilon$ we have $\tau_{B(x,h_d(\varepsilon))}(W_{T_{(x,\varepsilon)}})<\zeta_{T_{(x,\varepsilon)}}$ $\N_0$-a.e.

\noindent
We now bound $J_{21}$ using \reff{eq:N11}.
\begin{align*}
   J_{21}
&= 4\int_{A_\varepsilon} dx\;  \N_0\left[T_{(x,\varepsilon)}<\infty ; \int_{\tau_{B(x,h_d(\varepsilon))}}^{\zeta_{{T_{(x,\varepsilon)}}}}dt\; \N_{W_{T_{(x,\varepsilon)}}(t)}\left[T_{(x,\varepsilon)}<\infty ; \int_{0}^{L_{(x,\varepsilon)}}ds\;\ind_A (\hat W_s)\right]\right]\\
&\leq  4\int_{A_\varepsilon} dx\;  \N_0\left[T_{(x,\varepsilon)}<\infty ; \int_{\tau_{B(x,h_d(\varepsilon))}}^{\zeta_{{T_{(x,\varepsilon)}}}}dt\; \N_{W_{T_{(x,\varepsilon)}}(t)}\left[\int_{0}^{\sigma}ds\;\ind_A (\hat W_s)\right]\right]\\
&=4 \int_{A_\varepsilon} dx\;  \N_0\left[T_{(x,\varepsilon)}<\infty ; \int_{\tau_{B(x,h_d(\varepsilon))}}^{\zeta_{{T_{(x,\varepsilon)}}}}dt\; \int_A dy\; G({W_{T_{(x,\varepsilon)}}}(t),y)\right].
\end{align*}
Now we use \reff{eq:Nbeta}, the Cauchy-Schwarz inequality and \reff{eq:fk} to get
\begin{align*}
   J_{21}
&\leq 4 \int_{A_\varepsilon} dx\int_0^\infty  dt \;  \E_0\left[\tau_{B(x,\varepsilon)}>t\geq \tau_{B(x,h_d(\varepsilon))}; \int_A dy\; G(\beta_t,y)u_\varepsilon(\beta_t-x)\expp{-2\int_0^t u_\varepsilon(\beta_r-x)dr}\right]\\
&\leq 4 \left[2^{-1}\int_{A_\varepsilon} dx\;u_\varepsilon(x) \right]^{1/2} \left[\int_{A_\varepsilon} dx\int_0^\infty  dt \;  \E_0\left[t\geq \tau_{B(x,h_d(\varepsilon))}; \left(\int_A dy\; G(\beta_t,y)\right)^2\right]\right]^{1/2}\\
&\leq \cst \varphi_d(\varepsilon)^{-1/2} \left[\int_{A_\varepsilon} dx\;\P_0\left[\tau_{B(x,h_d(\varepsilon))}<\infty\right ] \sup_{x'\in B(0,2R_0)} \int dz\; G(z,x') \left(\int_A dy\; G(z,y)\right)^2\right]^{1/2}\\   
&\leq \cst \varphi_d(\varepsilon)^{-1/2} h_d(\varepsilon)^{(d-2)/2}.
\end{align*}
We used the strong Markov property at time $\tau_{B(x,h_d(\varepsilon))}$ and \reff{eq:intta} for the last two  inequalities. This implies that $J_{21}\leq \cste \varphi_d(\varepsilon)^{-1}h_d(\varepsilon)^{1/2}$.

\setcounter{cstejji}{\thecste}

Using the time reversal property of the Brownian snake, the strong Markov property at time $T_{(x,\varepsilon)}$  and \reff{eq:Egreen} we get
\begin{align*}
   J_{22}
&=4\int_{A_\varepsilon} dx\;  \N_0\Bigg[T_{(x,\varepsilon)}<\infty ; \int_0^{\tau_{B(x,h_d(\varepsilon))}}dt\; \N_{W_{T_{(x,\varepsilon)}}(t)}\left[T_{(x,\varepsilon)}<\infty ; \int_{T_{(x,\varepsilon)}}^{\sigma}ds\;\ind_A (\hat W_s)\right]\Bigg]\\
%&=4\int_{A_\varepsilon} dx\;  \N_0\Bigg[T_{(x,\varepsilon)}<\infty ; \int_0^{\tau_{B(x,h_d(\varepsilon))}}dt\; \\
%&\hspace{3cm} \N_{W_{T_{(x,\varepsilon)}}(t)}\left[T_{(x,\varepsilon)}<\infty ; \ce^*_{W_{T_{(x,\varepsilon)}}} \left[\int_{0}^{\sigma}ds\;\ind_A (\hat W_s)\right]\right]\Bigg]\\
&=8\int_{A_\varepsilon} dx\;  \N_0\Bigg[T_{(x,\varepsilon)}<\infty ; \int_0^{\tau_{B(x,h_d(\varepsilon))}}dt\; \\
&\hspace{3cm}\N_{W_{T_{(x,\varepsilon)}}(t)}\left[T_{(x,\varepsilon)}<\infty ; \int_0^{\zeta_{T_{(x,\varepsilon)}}}ds\int_{A} dy\; G(W_{T_{(x,\varepsilon)}}(s),y)\right]\Bigg].
\end{align*}
We will distinguish according to $\left\{s\geq \tau_{B(x,h_d(\varepsilon))}\right\}$ (integral $J_{23}$) and $\left\{s< \tau_{B(x,h_d(\varepsilon))}\right\}$ (integral $J_{24}$). We now bound $J_{23}$. We have
\begin{align*}
   J_{23}
&=8\int_{A_\varepsilon} dx\;  \N_0\Bigg[T_{(x,\varepsilon)}<\infty ; \int_0^{\tau_{B(x,h_d(\varepsilon))}}dt\; \\
&\hspace{2cm}\N_{W_{T_{(x,\varepsilon)}}(t)}\left[T_{(x,\varepsilon)}<\infty ; \int_{\tau_{B(x,h_d(\varepsilon))}}^{\zeta_{T_{(x,\varepsilon)}}}ds\int_{A} dy\; G(W_{T_{(x,\varepsilon)}}(s),y)\right]\Bigg]\\
&=8\int_{A_\varepsilon} dx \int_0^\infty  dt\;\E_0\Bigg[\tau_{B(x,h_d(\varepsilon))}> t; u_\varepsilon(\beta_t-x)\expp{-2\int_0^t u_\varepsilon(\beta_r-x)dr}\int_0^\infty ds\;\\
&\hspace{2cm}  \E_{\beta_t}\left[\tau_{B(x,\varepsilon)}>s\geq \tau_{B(x,h_d(\varepsilon))}; \int_{A} dy\; G(\beta_s,y) u_\varepsilon(\beta_s-x)\expp{-2\int_0^s u_\varepsilon(\beta_v-x)dv}\right]\Bigg]\\
&\leq \cst \varphi_d(\varepsilon)^{-1}\int_{A_\varepsilon} dx  \int_0^\infty  dt\;\E_0\Bigg[\tau_{B(x,h_d(\varepsilon))}>t; \val{\beta_t-x}^{2-d}\\
&\hspace{2cm}\left[\int_0^\infty ds\; \E_{\beta_t}\left[\tau_{B(x,\varepsilon)}>s; u_\varepsilon(\beta_s-x)^2\expp{-4\int_0^s u_\varepsilon(\beta_v-x)dv}\right]\right]^{1/2}\\
&\hspace{2cm}\left[\int_0^\infty ds\; \E_{\beta_t}\left[s\geq \tau_{B(x,h_d(\varepsilon))}; \left(\int_{A} dy\; G(\beta_s,y)\right)^2\right]\right]^{1/2}\Bigg]\\
&\leq \cst \varphi_d(\varepsilon)^{-1}\int_{A_\varepsilon} dx\;  \int_0^\infty  dt\;\E_0\Bigg[\tau_{B(x,h_d(\varepsilon))}> t; \val{\beta_t-x}^{2-d}\left[2^{-1}u_\varepsilon(\beta_t-x)\right]^{1/2}\\
&\hspace{2cm}\left[\E_{\beta_t}\left[\tau_{B(x,h_d(\varepsilon))}<\infty ; \E_{\beta_{\tau_{B(x,h_d(\varepsilon))}}}\left[\int_0^\infty ds\; \left(\int_{A} dy\; G(\beta_s,y)\right)^2\right]\right]\right]^{1/2}\Bigg]\\
&\leq \cst \varphi_d(\varepsilon)^{-3/2}\int_{A_\varepsilon} dx\;  \int_0^\infty  dt\;\E_0\left[\tau_{B(x,h_d(\varepsilon))}> t; \val{\beta_t-x}^{(6-3d)/2}
\P_{\beta_t}\left[\tau_{B(x,h_d(\varepsilon))}<\infty \right]^{1/2}\right]\\
&\hspace{2cm}\left[\sup_{x'\in B(0,2R_0)}\int dz'G(x',z')\left(\int_{A} dy\; G(z',y)\right)^2\right]^{1/2}\\
&\leq \cst \varphi_d(\varepsilon)^{-3/2}\int_{A_\varepsilon} dx\;  \int_{\val{z-x}\geq  h_d(\varepsilon)} dz\; G(0,z)\val{z-x}^{(6-3d)/2}h_d(\varepsilon)^{(d-2)/2}\val{z-x}^{(2-d)/2}.
\end{align*}
We used \reff{eq:Nbeta} twice for the second equality, \reff{eq:umajo} and Cauchy-Schwarz inequality for the first inequality, \reff{eq:fk} and the strong Markov property at time $\tau_{B(x,h_d(\varepsilon))}$ for the second and \reff{eq:intta} for the last. We easily deduce that $J_{23}\leq \cste \varphi_d(\varepsilon)^{-1}h_d(\varepsilon)$.
\setcounter{cstejjii}{\thecste}

\noindent
For $J_{24}$ we have using \reff{eq:Nbeta} twice and \reff{eq:umajo} twice,
\begin{align*}
  J_{24}
&=8\int_{A_\varepsilon} dx\;  \N_0\Bigg[T_{(x,\varepsilon)}<\infty ; \int_0^{\tau_{B(x,h_d(\varepsilon))}}dt\; \\
&\hspace{4cm}\N_{W_{T_{(x,\varepsilon)}}(t)}\left[T_{(x,\varepsilon)}<\infty ; \int_0^{\tau_{B(x,h_d(\varepsilon))}}ds\int_{A} dy\; G(W_{T_{(x,\varepsilon)}}(s),y)\right]\Bigg]\\
&=8\int_{A_\varepsilon} dx\;  \int_0^\infty  dt\int_0^\infty  ds\; \E_0\Bigg[\tau_{B(x,h_d(\varepsilon))}> t; u_\varepsilon(\beta_t-x)\expp{-2\int_0^t u_\varepsilon(\beta_r-x)dr}\\
&\hspace{4cm}\E_{\beta_t}\left[\tau_{B(x,h_d(\varepsilon))}>s; \int_{A} dy\; G(\beta_s,y) u_\varepsilon(\beta_s-x)\expp{-2\int_0^s u_\varepsilon(\beta_v-x)dv}\right]\Bigg]\\
&\leq \cst \varphi_d(\varepsilon)^{-2} \int_{A_\varepsilon} dx\;  \int_0^\infty  dt\int_0^\infty  ds\; \E_0\Bigg[\tau_{B(x,h_d(\varepsilon))}> t; \val{\beta_t-x}^{2-d}\\
&\hspace{4cm}\E_{\beta_t}\left[\tau_{B(x,h_d(\varepsilon))}>s; \int_{A} dy\; G(\beta_s,y)\val{\beta_s-x}^{2-d}\right]\Bigg]\\ 
&\leq \cst \varphi_d(\varepsilon)^{-2} M(d,h_d(\varepsilon)).
\end{align*}
Using \reff{eq:Mde} we get $J_{24}\leq  \cste   \varphi_d(\varepsilon)^{-1}h_d(\varepsilon)^\theta$.
As a conclusion we get
\[
 J_2\leq \cstejji \varphi_d(\varepsilon)^{-1}h_d(\varepsilon)^{1/2}+\cstejjii \varphi_d(\varepsilon)^{-1}h_d(\varepsilon)+\mcste \varphi_d(\varepsilon)^{-1}h_d(\varepsilon)^{\theta}.
\]

\textbf{Conclusion on the lower bound on $J$.}

 By combining the previous results, we get for $d\geq 4$,
 \[
 \varphi_d(\varepsilon)J\geq J_0- \cstejii  h_d(\varepsilon)^{\theta} -\varphi_d(\varepsilon) J_2\geq
 J_0- \cste 
 h_d(\varepsilon)^{\theta}.
 \]

 \subsection{End of the proof of proposition \ref{th:erange}}

 We deduce from formula (\ref{eq:N12}), that
 \[
    J_0=  \cstap \N_0\left[\left[\int_0^\sigma
    \ind_A(\hat W_s) ds \right]^2\right],\quad \text{and} \quad
    I_0= {\cstap}^2 
    \N_0\left[\left[\int_0^\sigma \ind_A(\hat W_s) ds \right]^2\right].
 \]
 Thus we get from section \ref{sec:majoI}  and \ref{sec:sublbJ} that for $\varepsilon$ small enough 
 \[
 \N_0\left[\left[\varphi_d(\varepsilon) \val{\range^\varepsilon\cap
 A_\varepsilon} - \cstap  \int_0^\sigma
 ds \;\ind_A(\hat W_s) \right]^2\right]
 \leq 
 \csteiii  h_d(\varepsilon)^{\theta/2}+2\mcste h_d(\varepsilon)^{\theta}.
 \]
 Take $\kappa<\theta/2$ and $\varepsilon_0$ small  to
get the second  upper bound of  proposition \ref{th:erange}.  
\findemo

\section{Capacity equivalence for the support and the range of $X$}% super-Brownian motion}
\label{sec:cap}
\initcste
Let $f:(0,\infty )\rightarrow [0,\infty
)$ be a decreasing function. We  put $f(0)=\lim_{r\downarrow 0} f(r) \in [0,\infty ]$. We define
the energy of a Radon measure $\nu$ on $\R^d$ with respect to the
kernel $f$ by: $ \ci_f(\nu)=\iint f(\val{x-y}) \nu(dx)\nu(dy)$,
and the  capacity of a set $\Lambda\in \cb(\R^d)$ by 
\[
\capaf{\Lambda} =\left[\inf_{\nu(\Lambda)=1} \ci_f(\nu)\right] ^{-1}.
\]
Following \cite{pp:gwt}, we say that two sets $\Lambda_1$ and
$\Lambda_2$ are capacity-equivalent 
if there exist two positive constants $c$ and $C$ such that for every kernel $f$, we have
\[
c\capaf{\Lambda_1} \leq \capaf{\Lambda_2}\leq C\capaf{\Lambda_1} .
\] 
The next lemma is an
immediate 
consequence of  the remarks in \cite{pps:tsbm} p.385.
\begin{lem}
\label{rem:majocapf}
Let $\Lambda\subset \R^d$ be a bounded Borel set. Suppose  there
exist two positive constants $c'$ and $\gamma$
such that 
\begin{equation*}
  \lim_{\varepsilon\rightarrow 0} \varepsilon^{\gamma-d}\val{\Lambda^\varepsilon}=c'.
\end{equation*}
Then there exists a constant $C$ such that for every kernel $f$, we have
\[
\capaf{\Lambda} \leq C\left[\int_0^1 f(r) r^{\gamma-1} dr\right]^{-1}.
\]  
\end{lem}
\noindent
For every measure $\mu\in M_f$, we set
\[
S_\varepsilon(\mu)=\iint \mu(dx)\mu(dy)\; p(\varepsilon^2,x-y),
\]
where $p$ is the Brownian transition density in $\R^d$:  
$p(t,x)=(2 \pi t)^{-d/2} \expp{-\val{x}^2/2t}$, $(t,x)\in (0,\infty )\times \R^d$. 
The next lemma is also an immediate consequence of   \cite{pps:tsbm} (p.387).
\begin{lem}
\label{rem:minocapf}
Let $\Lambda\subset \R^d$ be a bounded Borel set. Suppose there exist two  positive constants $c'$ and $\gamma$
and a measure $\mu\in M_f$ 
such that $\mu(\Lambda^c)=0$ and 
\[
  \lim_{\varepsilon\rightarrow 0} \varepsilon^{d-\gamma}
  S_\varepsilon(\mu)=c'.
\]
Then there exists a constant $c$ such that for every kernel $f$, we have
\[
c\left[\int_0^1 f(r) r^{\gamma-1} dr\right]^{-1}\leq \capaf \Lambda .
\] 
\end{lem}
For example, for every integer  $p\leq d$, we can consider the cube
$[0,1]^p$ as a subset of $\R^d$, and then we obviously have 
\[
\lim_{\varepsilon\rightarrow 0}
\varepsilon^{p-d}\val{\left([0,1]^p\right)^\varepsilon}=2 \pi^{(d-p)/2}/\Gamma((d-p)/2),
\]
and if $\mu$ is  Lebesgue measure on $[0,1]^p$,
\[
\lim_{\varepsilon\rightarrow
  0}\varepsilon^{d-p}S_\varepsilon(\mu)=(2\pi)^{(p-d)/2}.
\]
Thus we deduce from lemma \ref{rem:majocapf} and \ref{rem:minocapf}
that there exist two positive constants $c'_p$, $C'_p$, such that for
every kernel $f$,
\begin{equation}
\label{eq:capr}
c'_p \left[\int_0^1 f(r) r^{p-1} dr\right]^{-1}\leq \capaf{[0,1]^p} \leq
C'_p \left[\int_0^1 f(r) r^{p-1} dr\right]^{-1}.
\end{equation}
We shall prove the following result on  super-Brownian motion and ISE.
\begin{prop}
\label{prop:capaequivX}
   \begin{itemize}
      \item[(i)] Assume $d\geq 3$. Let $t>0$, $\nu\in M_f$.
   $\P^X_\nu$-a.s. on $\left\{X_t\neq 0\right\}$, the set $\supp X_t$ is
   capacity-equivalent to $[0,1]^2$.
      \item[(ii)] Assume $d\geq 5$. Let $t>0$, $\nu\in M_f$.
   $\P^X_\nu$-a.s. on $\left\{X_t\neq 0\right\}$,  the set $\range_t(X)$ is
   capacity-equivalent to $[0,1]^4$.
   Furthermore, if there exists a positive number
   $\rho<4$ such that $\lim_{\varepsilon\rightarrow 0}
   \varepsilon^{\rho-d} \val{(\supp \nu)^\varepsilon}=0$, then
   $\P^X_\nu$-a.s.  the set $\range_0(X)$ is 
   capacity-equivalent to $[0,1]^4$.
      \item[(iii)] Assume $d\geq 5$. The set $\range_t(W)$ is
   capacity-equivalent to $[0,1]^4$ $\N^{(1)}_0$-a.s.
   \end{itemize}
\end{prop}

\noindent
\textbf{Proof} of proposition \ref{prop:capaequivX} (i). Let $d\geq
3$. It is  well-known  that for $t>0$, $\P^X_\nu$-a.s. the set $\supp
X_t$ is 
bounded. Thus, thanks to theorem \ref{th:tribe}, $\P^X_\nu$-a.s., we have
\[
\lim_{\varepsilon\rightarrow 0} \varepsilon^{2-d}\val{\left(\supp
X_t\right)^\varepsilon} =\cstT (X_t,\ind).
\]
Now  apply lemma \ref{rem:majocapf} to $\Lambda=\supp X_t$, with
$\gamma=2$ and take $p=2$ in
\reff{eq:capr}. We get that
$\P^X_\nu$-a.s., on $\left\{X_t\neq 0\right\}$, there exists a (random) 
constant $C_1>0$, such that for every kernel $f$,
\[
\capaf{\supp X_t}\leq C_1 \capaf{[0,1]^2}.
\]
For the second part of (i), we use  lemma \ref{lem:cvSYt} below. Recall
 notation $Y_t$ from  section \ref{sec:snake}. 
\begin{lem}
\label{lem:cvSYt}
   Fix $t>0$ and $x\in \R^d$, $d\geq 3$. Then we have 
\[
\lim_{\varepsilon\rightarrow 0} \varepsilon^{d-2} (2\pi)^{d/2}S_\varepsilon
(Y_t)=\frac{4}{d-2} (Y_t,\ind),
\]
where the convergence holds  $\N_x$-a.e. and in $L^2(\N_x)$.
\end{lem}
\noindent
Let us explain how the proof is completed using lemma \ref{lem:cvSYt}. Thanks to lemma  \ref{rem:minocapf}, the above lemma and (\ref{eq:capr}) imply that
$\N_x$-a.e. on $\left\{Y_t\neq 0\right\}$, there exists a positive
constant $c_1$ such that for every kernel $f$,
\[
\capaf{\supp Y_t} \geq c_1 \capaf{[0,1]^2}.
\]
Now remember that for $t>0$, under $\P^X_\nu$,   we can
write $X_t=\sum_{i\in I}Y_t(W^i)$, where
$\sum_{i\in I} \delta_{W^i}$ is a Poisson measure on $C(\R^+,\cw)$ with
intensity $\int \nu(dx)\; \N_x[\cdot]$. 
%Furthermore there is only a finite
%number of non zero term in the sum $\sum_{i\in I}Y_t(W^i)$.
On $\left\{X_t\neq 0\right\}$, there exists $i_0$ such that
$Y_t(W^{i_0})\neq 0$. Then we have $\supp Y_t(W^{i_0}) \subset \supp
X_t$. 
Thus the previous lemma entails that there exists a.s. a positive
constant $c_1(W^{i_0})$ such that for every kernel $f$,
\[
\capaf{\supp X_t}\geq \capaf{\supp Y_t(W^{i_0})}\geq c_1(W^{i_0})
\capaf{[0,1]^2}
\]
This completes the proof of (i).
\findemo

\noindent
\textbf{Proof} of proposition \ref{prop:capaequivX} (ii). Let $d\geq 5$. We argue as in the proof of (i) using theorem  \ref{th:rangeX} instead of theorem \ref{th:tribe} and the following lemma instead of lemma \ref{lem:cvSYt}.

\begin{lem}
\label{lem:cvSrYt}
   Fix $t\geq 0$ and $x\in \R^d$, $d\geq 5$. Then we have for every $T>t\geq 0$,
\[
\lim_{\varepsilon\rightarrow 0} \varepsilon^{d-4} (2\pi)^{d/2} S_\varepsilon
\left(\int_t^T ds\;Y_s  \right) =\frac{16}{(d-2)(d-4)} \int_t^T ds \;(Y_s ,\ind),
\]
where the convergence holds $\N_x$-a.e. and in $L^2(\N_x)$.
\end{lem}
\findemo

\noindent
\textbf{Proof} of proposition \ref{prop:capaequivX} (iii).  Let $d\geq 5$. For the first part we argue as in the proof of (i) using the second part of corollary \ref{cor:cvpsrangeY} instead of theorem \ref{th:tribe}. Notice  that thanks to \reff{eq:Nnormal} and the scaling property of the family $(\N^{(r)}_0,r>0)$, the convergence in lemma \ref{lem:cvSrYt} also holds $\N^{(1)}_0$-a.s. The second part of (iii) is then a direct consequence of lemma \ref{rem:minocapf} (with $\mu=\int_0^T ds\; Y_s$ and $\gamma=4$) and \reff{eq:capr} (with $p=4$).
\findemo

\noindent
The proofs  of lemma \ref{lem:cvSYt} and lemma \ref{lem:cvSrYt} are
very similar. We shall only prove the latter. The former uses the same
techniques in a simpler way.

\noindent
\textbf{Proof} of lemma \ref{lem:cvSrYt}. We first want to show  the
convergence in $L^2(\N_x)$. 
Fix $T>t\geq 0$. 
By standard monotone class arguments, we deduce from the results of
section \ref{sec:annexemo}  an explicit expression for
\[
\N_x\left[\int_0^T \cdots \int_0^T ds_1\ldots ds_4\; \int\cdots \int
  Y_{s_1}(dx_1)\ldots Y_{s_4}(dx_4)
  g(s_1,\ldots,s_4,x_1,\ldots,x_4)\right], 
\]
where $g$ is any  measurable positive function on $(\R^+)^4\times (\R^d)^4$. Specializing to the case
$g(s_1,\ldots,s_4,x_1,\ldots,x_4)=\prod_{i=1}^4 \ind_{[t,T]}(s_i)
p(\varepsilon^2, x_1-x_2) p(\varepsilon^2, x_3-x_4)$, we get
\begin{multline*}
   \N_x \left[S_\varepsilon\left(\int_t^T ds\;Y_s\right)^2\right]\\
   \begin{aligned}[b]
     &= \inv{3} 4!2^3 \int_0^T ds \int dy \;p(s,x-y) \Bigg\{ 4\int^{T-s}_{(t-s)_+}
     ds_1\int dy_1\; p(s_1,y-y_1)\int_0^{T-s}ds_2\\
     &\hspace{1.5cm}\int dy_2\;
     p(s_2,y-y_2)\int_{(t-s-s_2)_+}^{T-s-s_2} ds_3\int dy_3 \;
     p(s_3,y_2-y_3)\int_0^{T-s-s_2} ds_4\\
     &\hspace{1.5cm}\int dy_4 \; p(s_4,y_2-y_4)
     %\\
     %&\hspace{1.5cm}
     \int_{(t-s-s_2-s_4)_+}^{T-s-s_2-s_4} ds_5\int dy_5 \; p(s_5,y_4-y_5)\\
     &\hspace{1.5cm}\int_{(t-s-s_2-s_4)_+}^{T-s-s_2-s_4} ds_6\int dy_6 \; p(s_6,y_4-y_6) \\
     &\hspace{3cm}[p(\varepsilon^2,y_1-y_3)p(\varepsilon^2,y_5-y_6)+
     p(\varepsilon^2,y_1-y_5)p(\varepsilon^2,y_3-y_6)\\
     &\hspace{7cm}+p(\varepsilon^2,y_1-y_6)p(\varepsilon^2,y_3-y_5)]  \\
     &\hspace{1cm}+ \int_0^{T-s}ds_7\int dy_7\;
     p(s_7,y-y_7)\int_{(t-s-s_7)_+}^{T-s-s_7}ds_8\int dy_8\; 
     p(s_8,y_7-y_8)\\
     &\hspace{1.cm}\phantom{+ \int_0^{T-s}ds_7\int dy_7\;
     p(s_7,y-y_7)} \int_{(t-s-s_7)_+}^{T-s-s_7} ds_9\int dy_9 \; p(s_9,y_7-y_9)\\
     &\hspace{1.5cm}\int_0^{T-s    } ds_{10}\int dy_{10} \; p(s_{10},y-y_{10})
     \int_{(t-s-s_{10})_+}^{T-s-s_{10}} ds_{11}\int dy_{11} \; p(s_{11},y_{10}-y_{11})\\
     &\hspace{1.5cm} \phantom{\int_0^{t-s    } ds_{10}\int dy_{10} \;
     p(s_{10},y-y_{10})} \int_{(t-s-s_{10})_+}^{T-s-s_{10}}  ds_{12}\int dy_{12} \;
     p(s_{12},y_{10}-y_{12}) \\ 
   \end{aligned}\\ 
   [
   p(\varepsilon^2,y_{8}-y_{9})p(\varepsilon^2,y_{11}-y_{12})+
   p(\varepsilon^2,y_{8}-y_{11})p(\varepsilon^2,y_{ 9}-y_{12})\\
   +
   p(\varepsilon^2,y_{8}-y_{12})p(\varepsilon^2,y_{ 9}-y_{11})]  \Bigg\}
\end{multline*}
We write  $J_1$, $J_2$, $J_3$, $J_4$, $J_5$, and $J_6$,  respectively
for the 
integrals corresponding to the 
integrands $p(\varepsilon^2,y_1-y_3)p(\varepsilon^2,y_5-y_6)$,
$p(\varepsilon^2,y_1-y_5)p(\varepsilon^2,y_3-y_6)$,
$p(\varepsilon^2,y_1-y_6)p(\varepsilon^2,y_3-y_5)$,
$p(\varepsilon^2,y_{8}-y_{9})p(\varepsilon^2,y_{11}-y_{12})$,
$p(\varepsilon^2,y_{8}-y_{11})p(\varepsilon^2,y_{ 9}-y_{12})$, and 
$p(\varepsilon^2,y_{8}-y_{12})p(\varepsilon^2,y_{ 9}-y_{11})$  respectively. 
As we shall see the integral $J_4$ gives the main contribution.  Before
proceeding to the calculations, we give  three useful bounds:
for every positive real numbers  $s$, $\varepsilon^2<2^{-1}(T^{-1}\wedge
  T)$, we have  for $d\geq 5$ 
\begin{align}
   \label{eq:2-d} 
   \int_0^T \left(\varepsilon^2+s+r\right)^{-d/2}dr
   &\leq
   \frac{2}{d-2}\left(\varepsilon^2+s\right)^{1-d/2}, \\
   \label{eq:4-d} 
   \int_0^T \left(\varepsilon^2+s+r\right)^{1-d/2}dr
   &\leq
   \frac{2}{d-4}\left(\varepsilon^2+s\right)^{2-d/2}, \\
   \label{eq:6-d} 
   \int_0^T \left(\varepsilon^2+r\right)^{2-d/2}dr
   &\leq H_T(\varepsilon):=
   \begin{cases}
      {2}{(d-6)^{-1}}
      \varepsilon^{6-d} & \text{if $d\geq 7$},\\
       4\ln {\varepsilon}^{-1} & \text{if $d=6$},\\
      \displaystyle \sqrt{6T}  & \text{if $d=5$}.
   \end{cases}
\end{align}
From now on, we assume that $\varepsilon^2<2^{-1}(T^{-1}\wedge
  T)$ and also $\varepsilon^2\ln \varepsilon^{-1} <T$ if $d=6$. Let us derive an upper bound on $J_1$. By repeated applications of the
Chapman-Kolmogorov identities, we get 
\initre
\begin{align*}\initer
   J_1
   & \leq  2^8 \int_0^T \dotsi \int_0^T ds\dots ds_6  \int dy \;p(s,x-y)\int dy_1\; p(s_1,y-y_1)
   \\
   &\hspace{1cm} \int dy_2\;
   p(s_2,y-y_2)\int dy_3 \; p(s_3,y_2-y_3)\int dy_4 \;
   p(s_4,y_2-y_4)\\ 
   &\hspace{1.5cm}\int dy_5 \; p(s_5,y_4-y_5) \int dy_6 \;
   p(s_6,y_4-y_6) p(\varepsilon^2,y_1-y_3)p(\varepsilon^2,y_5-y_6)\\
   & \leq  2^8 \int_0^T \dotsi \int_0^T ds\dots ds_6\;
   p(\varepsilon^2+s_1+s_2+s_3,0)p(\varepsilon^2+s_5+s_6,0).
\end{align*}
We can apply
   (\ref{eq:2-d}), (\ref{eq:4-d}) and (\ref{eq:6-d})  to get:
\initre
\begin{align*}\initer
   J_1
   & \leq  \frac{2^8}{(2\pi)^d }  T \int_0^{T}ds_1
   \frac{4}{(d-2)(d-4)} (\varepsilon^2+s_1)^{2-d/2}\; \frac{4}{(d-2)(d-4)}
   T\varepsilon^{4-d}   \\
   & \leq  \cste
   T^2\varepsilon^{4-d} H_T(\varepsilon) ,
\end{align*}
where the constant $\mcste$ depends only on $d$. We can use the same method for $J_2$:
\setcounter{cstei}{\thecste}
\begin{align*}
   J_2
   & \leq  2^8 \int_0^T \dotsi \int_0^T ds\dots ds_6  \int dy \;p(s,x-y)\int dy_1\; p(s_1,y-y_1)
   \\
   &\hspace{1cm} \int dy_2\;
   p(s_2,y-y_2)\int dy_3 \; p(s_3,y_2-y_3)\int dy_4 \;
   p(s_4,y_2-y_4)\\ 
   &\hspace{1.5cm}\int dy_5 \; p(s_5,y_4-y_5) \int dy_6 \;
   p(s_6,y_4-y_6) p(\varepsilon^2,y_1-y_5)p(\varepsilon^2,y_3-y_6)\\
   & \leq  2^8 \int_0^T \dotsi \int_0^T ds\dots ds_6
   \int dz\;
   p(s_4,z)
   p(\varepsilon^2+s_1+s_2+s_5,z)p(\varepsilon^2+s_3+s_6,z),
\end{align*}
where we made the change of variables  $z=y_2-y_4$. Since
$p(\varepsilon^2+s_3+s_6,z)\leq p(\varepsilon^2+s_3+s_6,0)$ and
$p(\varepsilon^2+s_1+s_2+s_5,z)\leq p(\varepsilon^2+s_1+s_2+s_5,0)$ , we can
 argue as for $J_1$ to get:
\initre
\begin{align*}\initer
   J_2
   & \leq  2^8 \int_0^T \dotsi \int_0^T ds\dots ds_6\;
   p(\varepsilon^2+s_1+s_2+s_5,0)p(\varepsilon^2+s_3+s_6,0).\\   
   & \leq  \cstei
   T^2\varepsilon^{4-d} H_T(\varepsilon) .
\end{align*}
By symmetry, we get $J_2=J_3$. We want now to find an upper bound on
$J_4$. Using (\ref{eq:2-d}), (\ref{eq:4-d}) and (\ref{eq:6-d})  we  get:
\initre
\begin{align*}\initer
   J_4
   &=2^6\int_0^T ds \int dy \;p(s,x-y)
   \Bigg[\int_0^{T-s} ds_7\int_{(t-s-s_7)_+}^{T-s-s_7}
   ds_8 \int_{(t-s-s_7)_+}^{T-s-s_7}
   ds_9 \\
   &\hspace{1cm}\int dy_7\;
   p(s_7,y-y_7) \int dy_8\; 
   p(s_8,y_7-y_8) \int dy_9 \; p(s_9,y_7-y_9) p(\varepsilon^2, y_8-y_9)\Bigg]^2\\
   &=2^6\int_0^T ds 
   \Bigg[\int_0^{T-s} ds_7\int_{(t-s-s_7)_+}^{T-s-s_7}
   ds_8 \int_{(t-s-s_7)_+}^{T-s-s_7}
   ds_9 \;p(\varepsilon^2+s_8+s_9,0)\Bigg]^2\\   
   &\leq 2^6(2\pi)^{-d}\int_0^T ds 
   \Bigg[\int_0^{T-s} ds_7 \frac{4}{(d-2)(d-4)} \left[\varepsilon^2+
   2(t-s-s_7)_+\right]^{2-d/2} \Bigg]^2\\
   &= \frac{2^{10}}{(2\pi)^{d}\left[(d-2)(d-4)\right]^2} \\
   &\hspace{1cm} \int_0^T ds 
   \Bigg[\varepsilon^{4-d} [(T-s)-(t-s)_+]+\int_0^{(t-s)_+} ds_7  \left[\varepsilon^2+
   2(t-s-s_7)_+\right]^{2-d/2} \Bigg]^2\\
   &\leq \frac{2^{10}}{(2\pi)^{d}\left[(d-2)(d-4)\right]^2}\int_0^T ds 
   \Bigg[\varepsilon^{4-d} [(T-s)\wedge(T-t)]+2^{-1}H_{2T}(\varepsilon)
   \Bigg]^2\\
   &\leq \frac{2^{10}}{(2\pi)^{d}}\left[\frac{\varepsilon^{4-d}}{(d-2)(d-4)}\right]^2 
   \left[\frac{(T-t)^3}{3}   
   + (T-t)^2t \right] +\cste T^2\varepsilon^{4-d} H_T(\varepsilon),     
\end{align*}
where the constant $\mcste$ depends only on $d$.
We now compute an upper bound on $J_5$:
\begin{align*}
   J_5
   &\leq 2^6 \int_0^T\dotsi\int_0^T ds\dots ds_{12}  \int dy \;p(s,x-y) \int dy_7\;
   p(s_7,y-y_7) \\
   &\hspace{1cm} \int dy_8\; 
   p(s_8,y_7-y_8) \int dy_9 \; p(s_9,y_7-y_9) \int dy_{10} \;
   p(s_{10},y-y_{10}) \\
   &\hspace{2cm}\int dy_{11} \; p(s_{11},y_{10}-y_{11}) \int dy_{12} \;
   p(s_{12},y_{10}-y_{12})  
   p(\varepsilon^2,y_{8}-y_{11})p(\varepsilon^2,y_{9}-y_{12})\\
   &\leq  2^6 \int_0^T\dotsi\int_0^T ds\dots  ds_{12}
   %&\\ \hspace{1cm} 
   \int dz \;p(s_7+s_{10},z) p(\varepsilon^2+s_8+s_{11},z)
   p(\varepsilon^2+s_{9}+s_{12},z),
\end{align*}
where we made the change of variables  $z=y_{10}-y_7$. 
Since $p(\varepsilon^2+s_{9}+s_{12},z)\leq
p(\varepsilon^2+s_{9}+s_{12},0)$, and
$p(\varepsilon^2+s_7+s_8+s_{10}+s_{11},0)\leq
p(\varepsilon^2+s_7+s_8+s_{10},0)$, we can argue as for $J_1$, and get:
\[
   J_5
   \leq \cstei T^2 \varepsilon^{4-d}  H_T(\varepsilon).
\]
By symmetry we get $J_6=J_5$. Combining the previous bounds leads to 
\initre
%\begin{multline*}
\begin{equation*}
\initer
 %  \label{eq:majoY4}
   \N_x\left[S_\varepsilon \left(\int_0^t  ds\; Y_s  \right)^2\right]%\\
   \leq
   \frac{2^{10}}{(2\pi)^d}
   \left[\frac{\varepsilon^{4-d}}{(d-2)(d-4)}
   \right]^2\left[\frac{(T-t)^3}{3} + (T-t)^2 t \right] +\cste T^2
   \varepsilon^{4-d} H_T(\varepsilon),
\end{equation*}
%\end{multline*}
where the constant $\mcste$  depends only on $d$.

We shall now find a lower bound for $\N_x \left[ S_\varepsilon(\int_t^T
ds\; Y_s) \int_t^T  ds\;(Y_s,\ind) \right]$. Using similar arguments as
in the beginning of the proof, we get 
\begin{align*}
   I:=&\N_x \left[ S_\varepsilon\left(\int_t^T
   ds\; Y_s\right) \int_t^T  ds\;(Y_s,\ind)\right]\\
   =&\inv{3}
   3!2^3 \int_0^T ds \int dy \; p(  s,x-y) \int_{(t-s)_+}^{T-s} ds_1\int
   dy_1\; p(s_1,y-y_1) \\
   &%\hspace{2cm}
   \hspace{1cm}\int_0^{T-s} ds_2\int
   dy_2\; p(s_2,y-y_2) \int_{(t-s-s_2)_+}^{T-s-s_2} ds_3 \int dy_3\; p(s_3,y_2-y_3)
    \int_{(t-s-s_2)_+}^{T-s-s_2} ds_4 \\
   &%\hspace{2cm}
   \hspace{2cm}\int dy_4\; p(s_4,y_2-y_4) 
   \left[p(\varepsilon^2,y_1-y_3)+p(\varepsilon^2,y_1-y_4)+p(\varepsilon^2,y_3-y_4)\right] .
\end{align*}
Since we are looking for a lower bound, we restrict our attention to the
term $p(\varepsilon^2, y_3-y_4)$. We get
\initre
\begin{align*}\initer
   I
   &\geq 2^4 \int_0^T ds \int_{(t-s)_+}^{T-s}  ds_1\int_0^{T-s} ds_2
   \int_{(t-s-s_2)_+}^{T-s-s_2} ds_3 \int_{(t-s-s_2)_+}^{T-s-s_2} ds_4\;
   p(\varepsilon^2+s_3+s_4,0)  \\
   &\geq  \frac{2^4}{(2\pi)^{d/2}}\frac{4}{(d-2)(d-4)} 
   \int_0^T ds \; \left[(T-s)\wedge(T-t) \right]\\
   &\hspace{2cm}\int_0^{T-s} ds_2\left[\left(\varepsilon^2 +2(t-s-s_2)_+
   \right)^{2-d/2} - 2\left(\varepsilon^2+(T-s-s_2)\right)^{2-d/2} \right]\\
   &\geq  \frac{2^6}{(2\pi)^{d/2}}\frac{1}{(d-2)(d-4)} 
   \int_0^T ds \; \left[(T-s)\wedge(T-t) \right]
   \left[\varepsilon^{4-d}(T-s-(t-s)_+) -2H_T(\varepsilon)\right]\\
   &\geq  \frac{2^6}{(2\pi)^{d/2}}\frac{\varepsilon^{4-d}}{(d-2)(d-4)} 
   \left[\frac{(T-t)^3}{3} + (T-t)^2 t \right]-\cste T^2 H_T(\varepsilon),
\end{align*}
where $\mcste$ depends only on $d$. 
Finally we deduce from section \ref{sec:annexemo}, with
$\varphi(s)=\ind_{[0,T-t]}(s)$,  that 
\[
\N_x \left[\left[\int_t^T ds \;(Y_s,\ind)
  \right]^2\right]=4\left[\frac{(T-t)^3}{3} + (T-t)^2 t \right]. 
\]
Combining the previous results, we get for $\varepsilon$ small enough 
\initre
\begin{align*}\initer
   \N_x\left[\left[\varepsilon^{d-4}(2\pi)^{d/2}S_\varepsilon\left(\int_t^T ds\;Y_s \right)
   -\frac{2^4}{(d-2)(d-4)} \int_t^T ds\;
   (Y_s,\ind)\right]^2\right]
   &\leq \cste T^2 \varepsilon^{d-4}H_T(\varepsilon)\\
   &\leq \cste T^2 \varepsilon,
\end{align*}
where $\mcste$ depends only on $d$. This gives  the  convergence in
$L^2(\N_x)$. 
Now $S_\varepsilon\left(\int_t^T ds\;Y_s \right)$ is monotone decreasing
in $\varepsilon$ (cf lemma 5.3 in \cite{pps:tsbm}). The 
$\N_x$-a.e. convergence then follows from the previous estimate by an
application of the Borel-Cantelli lemma and monotonicity arguments.
\findemo

\section{Appendix}
%\label{sec:annexe}
\subsection{Formula for moments of the Brownian snake}
\label{sec:annexemo}
For the reader's convenience, we recall some explicit formulas for
moments of the Brownian snake. These formulas are well-known, at least
in the context of superprocesses (see e.g. Dynkin
\cite{d:rfsmsi}). 
We can compute the Laplace functional of $\int_0^t ds (Y_s,\varphi(s))$
for $\varphi\in \cb_{b+}(\R^+\times \R^d)$. To this 
end start from  the
finite dimensional Laplace functional (\ref{eq:fLaplaceX}) with
$t_i=i/m$, $\varphi_i=\inv{m}\varphi(i/m)$ for a nonnegative continuous
function $\varphi$
with compact support  on $\R^+\times \R^d$. Thanks to
the continuity of the process $X$, by a suitable passage to the limit,
we get for $\nu\in M_f$
\[
\E_\nu^X\left[\exp{\left[-\int_0^t
(X_{t-s},\varphi(s))ds\right]}\right] =
\exp{\left[-(\nu,v(t))\right]} ,
\]
where $v$ is a nonnegative solution of (\ref{eq:integrale}) with
right-hand side $J(t,x)=\int_0^t ds\;P_{t-s}[\varphi(s)](x)$. This can be
extended by monotone class arguments  to any $\varphi\in \cb_{b+}(\R^+\times
\R^d)$. The uniqueness of the solution is easily established using
arguments similar to the classical Gronwall lemma.
Then we get $v(t,x)=
\N_x\left[1-\exp{\left[-\int_0^t ds\; (Y_{t-s},
      \varphi(s))\right]}\right]$, thanks to theorem \ref{th:XY}.

Now we
introduce an auxiliary power series. Let us consider the 
analytic function $f(\lambda)=1-\sqrt{1-\lambda}$ for
$\val{\lambda}<1$. It is easy to check that for $\val{\lambda}<1$, we
have 
\[
f(\lambda)=\sum_{n=1}^\infty  \gamma_n\lambda^n ,
\]
where the sequence $(\gamma_n,n\geq 1)$ is defined by $\gamma_1=1/2$ and
the recurrence 
\[
\gamma_n=\inv{2}\sum_{k=1}^{n-1} \gamma_k\gamma_{n-k}\quad \text{for}\quad n\geq 2
\]
(use the fact that $f$ solves $2f(\lambda) =f(\lambda)^2+\lambda$). Now let $T>0$
and $J$
a nonnegative measurable function on $\R^+\times \R^d$, 
such that $M_T=\sup_{[0,T]\times \R^d} J(t,x)<\infty $. We 
define the family of measurable functions $(h_n,n\geq 1)$ on $\R^+\times
\R^d$, by the initial condition 
\begin{align}
   \nonumber
   h_1(t)&=J(t),\\
\intertext{and the recurrence}
   \label{eq:rechn}
   h_n(t)&=2\sum_{k=1}^{n-1}\int_0^t ds\; P_s\left[h_k(t-s) h_{n-k}(t-s)\right]\quad \text{for} \quad n\geq 2 .
\end{align}
We clearly have for every $n\geq 1$,
\[
\sup_{[0,T]\times \R^d} \val{h_n}\leq [4T]^{n-1} [2M_T]^n \gamma_n.
\]
Thus the power series $w(\lambda,t)=\sum (-1)^{n+1}\lambda^n h_n(t)$ is normally convergent on
$[0,T]\times \R^d$ for $\val{\lambda}<[8TM_T]^{-1}$. And it clearly
solves the integral equation on $[0,T]\times \R^d$
\begin{equation}
   \label{eq:vlambda}
   w(t)+2\int_0^t ds\; P_s \left[w(t-s)^2\right]=\lambda
   J(t).
\end{equation}
To get the uniqueness  of the solution to the previous integral equation, use
arguments similar to Gronwall's lemma. Finally we can compute the moments
for the process $Y$ under $\N_x$. Indeed, let $\varphi\in \cb_{b+}
(\R^+\times \R^d)$. We have shown that for $\lambda>0$, the function
$v_\lambda(t,x)=\N_x \left[1-\exp{-\lambda \int_0^t (Y_{t-s},\varphi(s))}
  ds\right]$ is the unique solution to (\ref{eq:vlambda}) on $\R^+\times
\R^d$ with $J(t,x)=\int_0^t ds\; P_{s} [\varphi(t-s)](x)$. Thus for
$\lambda\geq 0$ small enough, we have 
$v_\lambda(t)=w(\lambda,t)$. Then from the series expansion for
$w(\lambda,t)$, we get for every integer $n\geq 1$
\[
\N_x\left[\left(\int_0^t ds\;
    (Y_{t-s},\varphi(s))\right)^n\right]= n!h_n(t,x),
\]
where the functions $h_n$ are defined by 
$h_1(t)=\int_0^t ds\; P_{s}[\varphi(t-s)]$, 
and the recurrence (\ref{eq:rechn}). In the same way it can be shown
that for every $\varphi\in \cb_{b+}(\R^d)$, for every $t\geq 0$, $n\geq 1$,
\[
\N_x\left[(Y_{t},\varphi)^n\right]= n!h_n(t,x),
\]
where the functions are defined by 
$h_1(t)=P_{t}[\varphi]$,
and the recurrence (\ref{eq:rechn}).

\subsection{Some properties of the function $u_1$}
\label{sec:annexeu}
We consider the function $u_1$,  which is the maximal solution on $(1,\infty )$ of
the non linear differential equation
\[
{u}''(r)+\frac{d-1}{r} u'  (r)=4u(r)^2.
\]
\begin{lem}
   \label{lem:a0}
There exist positive constants $\csta$, $\cstb$
and $\cstaaa$, depending only on $d$, such
that 
\begin{align}
\nonumber   
\lim_{r \rightarrow \infty }r^{d-2}u_1(r)=\csta\quad \text{if $d\geq 5$}, 
\quad 
& \lim_{r \rightarrow \infty }r^2\log(r)\;u_1(r) =\csta=1/2 \quad \text{if $d=4$};\\
\intertext{furthermore for every $r>1$,} 
\label{eq:umino1}
u_1(r) \geq \csta r^{2-d}\quad \text{if $d\geq 5$},
\quad 
& u_1(r) \geq \csta r^{-2} \log (2r)^{-1}\quad \text{if $d=4$};
\intertext{and for  every $r\geq 4/3$,} 
 \label{eq:umajo1}
   u_1(r)\leq \cstb r^{2-d}\quad \text{if $d\geq 5$},
\quad 
& u_1(r)\leq \cstb [2 r^2 \log (r)]^{-1}\quad \text{if $d=4$}; 
\end{align}
\begin{align}
\label{eq:udel1}
   u_1(r)&\leq \csta r^{2-d}+\cstaaa r^{6-2d}\quad\text{if $d\geq 5$},\\
\label{eq:udel41}
 u_1(r)&\leq \csta r^{-2} \log (r)^{-1}+\cstaaa r^{-2} \log (r)^{-2}\log(\log (r))\quad \text{if $d=4$}.
\end{align}
\end{lem}
We will prove this lemma by giving the asymptotic expansion of $u_1$ at $\infty $. For $d\geq 5$, we will see  the constant $\csta$ can be expressed as the radius of convergence of a series.
We introduce  the auxiliary function 
\[
z(t)=4(d-2)^{-2(d-1)/(d-2)} t\;
u_1\left[\left(\frac{t}{d-2}\right)^{1/(d-2)}\right],
\quad \text{for } t>d-2.
\]
This function is solution on $(d-2,\infty )$ of 
\begin{equation}
   \label{eq:diffy}
y''(t) =t^{-\delta-2} y(t)^2,
\end{equation}
where $\displaystyle \delta=\frac{d-4}{d-2}$. 
We deduce from \cite{t:abs} p.132 that the function $z(t)$ (i.e. $r^{d-2}u_1(r)$) is decreasing.

\noindent
\textbf{Proof} in the case $d\geq 5$. We deduce from theorem 1.1 of \cite{t:abs} that the limit $q=\lim_{t\rightarrow \infty } z(t)$ exists and is positive. Hence by integrating (\ref{eq:diffy}) twice from $t$ to
$\infty $, we get for $t>d-2$,
\begin{equation}
   \label{eq:integz}
z(t) -q =\int_t^\infty  (r-t) r^{-\delta-2} z(r)^2 dr.
\end{equation}
Now consider the sequence $(q_n,n\geq 0)$ defined by $q_0=1$ and the recurrence
\[
q_{n}=\inv{n\delta
  (n\delta+1)}\sum_{k=0}^{n-1} q_kq_{n-k-1} , \quad \text{for}\quad n\geq 1.
\]
Clearly we have for every $n\geq 0$,
$q_n\leq 2\left[4/\delta\right]^{n} \gamma_{n+1},$
where the sequence $(\gamma_n,n\geq 1)$ has been introduced in section
\ref{sec:annexemo}. The power series $\sum q_n q^{n+1} t^{-\delta n}$ is
convergent and even $C^\infty $ as a function of $t$ at least for
$t>t_1=\left[4q/\delta\right]^{1/\delta}$. This power series also solves
 (\ref{eq:integz}) for $t>t_1$.   The same arguments as in the proof of
 the Gronwall lemma
show that equation 
(\ref{eq:integz}) possesses a unique  solution bounded in a
neighborhood of infinity. Thus the function $z$ and the power series agree for $t>t_1$.
Since the function $z$ is analytic on
$(d-2,\infty )$ and since $q$ and the coefficients $q_n$ are positive, we deduce that the radius of convergence of the series $\sum q_n s^n$ is $q (d-2)^{-\delta}$ and that  for $t>d-2$
\[
z(t) =\sum_{n=0}^\infty  q_n q^{n+1} t^{-n\delta}.
\]
Thus we get with obvious notation for $r>1$, 
\begin{align*}
   u_1(r)
   &= 4^{-1} (d-2)^{d/(d-2)} r^{2-d}\sum_{n=0}^\infty  q_n q^{n+1}
   (d-2)^{-n(d-4)/(d-2)} r^{-n(d-4)}\\
   &= r^{2-d}\sum_{n=0}^\infty  \delu_n
   r^{-n(d-4)}.
\end{align*}
Since the function $r^{d-2} u_1(r)$ is decreasing we get
\reff{eq:umino1}. Since the real numbers $(\delu_n, n\geq 0)$ are positive, 
 \reff{eq:umajo1} and \reff{eq:udel1}  follow easily. Notice that $4 (d-2)^{-2} \csta$ is the radius of convergence of the series $\sum q_n s^n$.
\findemo

\noindent
\textbf{Proof} in the case \textbf{$d=4$}. 
We write $f(t)\thicksim g(t)$ at $0+$ when the real function $f$ and $g$ are positive or negative on $I=(0,0+\varepsilon)$ for some $\varepsilon>0$ and 
$\lim_{t\in I, t\rightarrow 0} f(t)/g(t)=1$. We also write $f(t)\thicksim g(t)$ at $\infty $ when $f(1/t)\thicksim g(1/t)$ at $0+$.
We know from  \cite{t:abs} p.133 that $z(t)\thicksim  \log (t)^{-1}$ at $\infty $. We deduce from \reff{eq:diffy} that $z'(t)$ is negative on $(2,\infty )$ and $z''(t)\thicksim  [t\log (t)]^{-2}$ at $\infty $. By integration, we get $z'(t)\thicksim  t^{-1}\log (t)^{-2}$ at $\infty $. We now consider the function $w(s)=z(\expp{s})$ which solves $w''-w'=w^2$ on $(\log 2,\infty )$. Notice that  the function $w$ is positive decreasing and $w'$ is negative. We also have $w(s)\thicksim s^{-1}$, $w'(s)\thicksim -s^{-2}$ and $w''(s)=o(s^{-2})$ at $\infty$. Thus the function defined on $(0,\infty )$ by
\[
p(w(s))=w'(s), \quad \text{for}\quad s\in (\log 2, \infty ),
\]
is well defined and even of class $C^1$, and $p'(w(s))=w''(s)/w'(s)$. Thus the function $p$ can be extended as a $C^1$ function on $[0,\infty )$ by setting $p(0)=0$ and $p'(0)=0$. Furthermore it solves
\[
p(w)p'(w)-p(w)=w^2 \quad \text{on}\quad [0,\infty ).
\]
We also have $p(w)\thicksim -w^2$ at $0+$. We consider the sequence $(\rho_n,n\geq 2)$ defined by $\rho_2=1$ and the recurrence 
\[
\rho_n=\sum_{k=2}^{n-1} k \rho_k\rho_{n-k+1},\quad \text{for} \quad n\geq 3.
\]
The radius of convergence of the series $\sum (-1)^{n+1} \rho_n w^n$ is $0$, nevertheless we will prove it is the asymptotic expansion of $p$ at $0+$. We set $H_n(w)=\sum_{k=2}^n (-1)^{k+1} \rho_k w^k$ for $n\geq 2$. We now prove by recurrence that $p(w)=H_n(w)+h_n(w)$, where $h_n(w)=o(w^n)$ at $0+$. This is true for $n=2$. Let us assume it is true at stage $n$. 
Let $g_{n,\alpha}(w)= (1-\alpha)(-1)^{n}\rho_{n+1} w^{n+1} -h_n(w)$. We easily have 
\begin{align*}
   g'_{n,\alpha}(w)p(w)+g_{n,\alpha}(w)[H'_n(w)-1]
&=\alpha(-1)^n \rho_{n+1}w^{n+1}+o(w^{n+1}),\\
&=(-1)^{n+1} \rho_{n+2}w^{n+2}+o(w^{n+2}), \quad \text{if} \quad \alpha=0.
\end{align*}
Let us assume $n$ is even. For $\alpha=0$, the above right hand side is negative on $(0,\varepsilon]$, for $\varepsilon$ small enough. Since $p$ is negative and $[H'_n(w)-1]<0$ on $[0,\varepsilon]$, for $\varepsilon$ small, we see that $g_{n,0}(w)<0$ implies $g_{n,0}'(w)\geq 0$. As $g_{n,0}(0)=0$, we get by contradiction that $g_{n,0}\geq 0$ on $[0,\varepsilon]$. This implies $h_n(w)\leq \rho_{n+1} w^{n+1}$. 
Similar arguments for $\alpha>0$ implies that $g_{n,\alpha}\leq 0$ on $[0,\varepsilon_\alpha]$ for $\varepsilon_\alpha>0$ small enough. Since this holds for any $\alpha>0$ and since $h_n(w)\leq \rho_{n+1} w^{n+1}$ for $w$ small enough,  we deduce that $h_{n+1}(w)=h_n(w)-\rho_{n+1} w^{n+1}=o(w^{n+1})$. If $n$ is odd the proof is similar.

From the definition of $p$, we then have $w'(s)=H_n(w(s))+O(w(s)^{n+1})$ at $\infty $. For $n=3$ this gives $w'(s)=-w(s)^2+2 w(s)^3+O(w(s)^4)$ at $\infty $. Since $w(s)\thicksim s^{-1}$  at $+\infty $, we deduce by integration that 
\[
\inv{w(s)}-2 \log w(s) + O(1)=s \text{ at infinity}.
\]
Standard arguments yields $w(s)=s^{-1}+2 s^{-2} \log (s) + O(s^{-2})$ at infinity. Thus we have
\[
u_1(r)=\inv{r^2} \left[\inv{2\log (r)} +\frac{\log (\log (r))}{4\log (r)^2} +O\left(\log (r)^{-2} \right)\right] \text{ at } +\infty .
\]
Notice the previous calculation can be continued to give an asymptotic expansion of $u_1$ at infinity. The inequalities
 \reff{eq:umajo1} and \reff{eq:udel41}  follow easily.
We will now prove that for every $r>1$,
$u_1(r) \geq [2 r^2 \log (2r)]^{-1}$. We consider the function $w(r)=u_1(r)-[2r^2 \log(2r)]^{-1}$. The function $w$ is positive at least over $(1,1+\eta)\cap(\eta^{-1}, \infty )$ for $\eta$ small. Let us assume that $w$ achieves its minimum at $r_0$ and that $w(r_0)\leq 0$. Then we have $r_0\in [1+\eta,\eta^{-1}]$, $w'(r_0)=0$ and $w''(r_0)\geq 0$. An easy computation gives 
\[
w''(r)=4w(r)\left[u_1(r)+\inv{2r^2\log (2r)}\right]-\frac{3}{r} w'(r)-\inv{2 r^4 (\log (2r))^3}.
\]
Evaluation at $r=r_0$ implies that $w''(r_0)<0$. This contradicts the assumption. Hence $w$ is positive, that is we get \reff{eq:umino1} for $d=4$.
\findemo

\end{document}